\documentclass[9pt,reqno]{amsart}
\usepackage{graphicx}
\usepackage{verbatim}
\usepackage{textcomp}
\usepackage{amssymb}
\usepackage{cite}
\usepackage{amsmath}
\usepackage{latexsym}
\usepackage{amscd}
\usepackage{amsthm}
\usepackage{mathrsfs}
\usepackage{bm}
\usepackage{url}
\usepackage{hyperref}
\usepackage{bookmark}
\allowdisplaybreaks[3]
\newtheorem{thm}{Theorem}[section]

\newtheorem{lem}[thm]{Lemma}
\newtheorem{prop}[thm]{Proposition}

\theoremstyle{definition}

\theoremstyle{remark}

\numberwithin{equation}{section}
\begin{document}
\title[A weighted Reilly type integral formula for differential forms]
{A weighted Reilly type integral formula for differential forms and its applications}

\author{Liyi Cao }
\author{Guangyue Huang}
\author{Hongru Song$^*$}

\address{School of Mathematics and Statistics, Henan Normal University, Xinxiang 453007, P.R.China  }

\email{lycao2025@126.com(L. Cao) }
\email{hgy@htu.edu.cn(G. Huang) }
\email{songhongru@htu.edu.cn(H. Song) }

\thanks{The research of the second author is supported by key projects of the Natural Science Foundation of Henan Province (No. 252300421303) and NSFC (No. 11971153), the third author is supported by Postdoctoral Fellowship Program of CPSF (No. GZC20230734).}
\thanks{$^*$Corresponding author: songhongru@htu.edu.cn(H. Song) }

\begin{abstract}
In this paper, we derive a weighted Reilly type integral formula for differential forms on a compact smooth metric measure space with boundary. As applications, a lower bound of the spectrum for the weighted Hodge Laplacian acting on differential forms on the boundary, and some special properties for $p$-th absolute cohomology space with respect to the lowest $p$-curvatures of the boundary have been obtained, respectively. Furthermore, we obtain a lower bound for the first positive eigenvalue of the Steklov eigenvalue problem on differential forms which is related to the lowest principal curvature of the boundary, and a comparison result between the eigenvalues of the Steklov eigenvalue problem and the
Hodge Laplacian on the boundary. On the other hand, for closed submanifolds of weighted Euclidean space, we derive universal inequalities for the sum of eigenvalues with respect to the weighted Hodge Laplacian, which can be seen as a generalization of Levitin-Parnovski inequality.
\end{abstract}



\maketitle

\section{Introduction}
Let $(M,g)$ be an $(n+1)$-dimensional compact Riemannian manifold with the smooth boundary $\partial M$, and $g$ denote the Riemannian metric, where $n\geq2$. For any differential $p$-form $\omega\in A^p(M)$, the Hodge Laplacian $\Delta^H$ is given by
\begin{equation}\label{2-Sec-1}
\Delta^H\omega=(d\delta+\delta d)\omega,
\end{equation}
where $d$ is the differential operator and $\delta=d^*$ is the codifferential operator, which is the adjoint of $d$ with respect to the canonical inner product of forms. In particular, if $\omega=u$ is a smooth function, then
\begin{equation}\label{2-Sec-111}
\Delta^Hu=-\sum_{A=1}^{n+1}u_{AA}:=\Delta u,
\end{equation}
where $\Delta$ is the Laplacian (which is also called the Laplace-Beltrami operator) on the smooth function $u$.

Let the triple $(M, g, d\mu)$ be a smooth metric measure space with the weighted measure
$$d\mu:=e^{-f}dv,$$
where $f$ is a fixed smooth real-valued function defined on $M$, and $dv$ denotes the volume form with respect to metric $g$.
The weighted Laplacian (which is also called the Witten Laplacian, see \cite{li05,lD11,HL2014}) associated with $f$
is defined by
\begin{equation}\label{2-Sec-222}
\Delta_fu:=-{\rm div}_f(\nabla u),
\end{equation}
where ${\rm div}_f(\nabla u)=e^f(e^{-f}u_A)_A=-\Delta u-\langle\nabla f,\nabla u\rangle$,
which is also adjoint with respect to the weighted measure $d\mu$.
That is,
$$\int_{M}u\Delta_f v\,d\mu=\int_{M}\langle\nabla u,\nabla v\rangle\,d\mu=\int_{M}v\Delta_f u\,d\mu, \ \ \ \forall \ u,v\in
C^\infty_0(M).$$
Associated with the weighted Laplacian, the $m$-dimensional Bakry-\'{E}mery Ricci curvature (see \cite{Bakry85,Bakry94,li05,Wei09,lW2015}) is
$${\rm Ric}_f^m={\rm Ric}+\nabla^2f-\frac{1}{m-n-1}df\otimes df,
$$ where $m\geq n+1$ is a constant, and $m=n+1$ if and only if $f$ is a constant. Define
\begin{equation}\label{2-Sec-333}
{\rm Ric}_f={\rm Ric}+\nabla^2f.
\end{equation}
Then ${\rm Ric}_f$ can be seen as the case of the $\infty$-dimensional Bakry-\'{E}mery Ricci curvature.

From \eqref{2-Sec-111}, we know that for any smooth function, the Hodge Laplacian is consistent with the Laplacian.
It is natural to consider the corresponding weighted Hodge Laplacian which is consistent with the Witten Laplacian given by \eqref{2-Sec-222}.

Let $J:\partial M\rightarrow M$ be the inclusion mapping and $J^*$ be the restriction of differential forms to the boundary. Let $N$ be the inner unit normal vector field along $\partial M$ and $S:\mathcal{X}(\partial M)\rightarrow \mathcal{X}(\partial M)$ be the Weingarten operator with respect to $N$, that is, for all $X\in\mathcal{X}(\partial M)$, we have $S(X)=-\nabla_XN$. Denote by $H=\frac{1}{n}{\rm tr}_g(S)$ the mean curvature of $\partial M$ and $da$ denotes the measure with respect to the induced metric on $\partial M$. For a (1,1) tensor $T$, its induced operator $T^{[p]}: A^p\rightarrow A^p$ is defined by
\begin{align}\label{2-Sec-5}
(T^{[p]}\omega)(X_1,\ldots,X_p):=\sum_{\alpha=1}^p\omega(X_1,\ldots,T(X_\alpha),\ldots,X_p).
\end{align}
Naturally, inspired by \eqref{2-Sec-222}, we can define the weighted codifferential operator $\delta_f$ associated with $f$ by
\begin{equation}\label{2-Sec-2}
\delta_f\omega=e^f\delta(e^{-f}\omega).
\end{equation}
Then, the weighted Hodge Laplacian $\Delta^H_f$ associated with $f$ is given by
\begin{equation}\label{2-Sec-3}
\Delta_f^H\omega=(d\delta_f+\delta_f d)\omega.
\end{equation}
By a direct calculation, we obtain the following Bochner formula (see Proposition \ref{1-sec-Prop-1}):
\begin{equation}\label{2-Sec-4}
\Delta_f^H=\nabla_f^*\nabla+W_f^{[p]},
\end{equation}
where $\nabla_f^*\nabla:=\nabla^*\nabla+{\nabla}_{\nabla f}$ and the $f$-Weitzenb\"{o}ck curvature $W_f^{[p]}$ is defined by $W_f^{[p]}:=W^{[p]}+(\nabla^2f)^{[p]}$.

In particular, by using \eqref{2-Sec-4}, we can establish the following weighted Reilly type formula for differential forms with the weighted measure:

\begin{thm}\label{1-Th-1}
Let $(M, g, d\mu)$ be a compact smooth metric measure space and $V$ be a smooth function defined on $M$. Then,
\begin{align}\label{1-Th-Formula-1}
\int_M&V(|\delta_f\omega|^2+|d\omega|^2-|\nabla\omega|^2)\,e^{-f}dv\notag\\
=&\int_M[-2\langle\omega,i_{\nabla V}(d\omega)\rangle +V\langle W_{f,V}^{[p]}(\omega),\omega\rangle]\,e^{-f}dv\notag\\
&+\int_{\partial M}[-V_N|J^*\omega|^2+2V\langle\delta^{\partial M}_f(J^*\omega),i_N\omega\rangle+ V\mathcal{B}_f(\omega,\omega)]\,e^{-f}da,
\end{align}
where $H_f=H+\frac{1}{n}f_N$, the $(f,V)$-Weitzenb\"{o}ck curvature $W_{f,V}^{[p]}$ is defined by
\begin{equation}\label{1-Th-Formula-2}
VW_{f,V}^{[p]}=VW_f^{[p]}+(\Delta_fV)g+(\nabla^2V)^{[p]}
\end{equation}
and
\begin{equation}\label{1-Th-Formula-3}
\mathcal{B}_f(\omega,\omega)=\langle S^{[p]}(J^*\omega),J^*\omega\rangle+nH_f|i_N\omega|^2-\langle S^{[p-1]}(i_N\omega),i_N\omega\rangle.
\end{equation}
In particular, \eqref{1-Th-Formula-3} is equivalent to
\begin{equation}\label{1-Th-Formula-4}
\mathcal{B}_f(\omega,\omega)=\langle S^{[p]}(J^*\omega),J^*\omega\rangle+\langle S^{[n+1-p]}(J^* *\omega),J^* *\omega\rangle+f_N|i_N\omega|^2.
\end{equation}
\end{thm}

When $f$ is a constant, our Theorem \ref{1-Th-1} becomes Xiong's Theorem 1 in \cite{Xiong2024}. Thus, our Theorem 1.1 can be seen as a generalization to Theorem 1 of Xiong in \cite{Xiong2024};
When $\omega=du$ and $f$ is a constant, where $u$ is a smooth function, then \eqref{1-Th-Formula-1} becomes
a special case of Qiu and Xia's Theorem 1.1 in \cite{QX2015}. Moreover, in this case, the $(f,V)$-Weitzenb\"{o}ck curvature $W_{f,V}^{[1]}$ becomes $\infty$-dimensional generalized Ricci curvature given by (1.3) in \cite{HMZ2024}, which is a generalization to (0,2)-tensor $Q$ defined in \cite{lX2019}; When $V=1$ and $f$ is a constant, our \eqref{1-Th-Formula-1} becomes the formula (4) of Raulot and Savo in \cite{Savo2011}. That is to say, for $V=1$, our Theorem \ref{1-Th-1} provides a weighted version to Theorem 3 in \cite{Savo2011} under the metric measure space.

Let $\eta_1(x),\ldots,\eta_n(x)$ be the principal curvatures of $\partial M$ at $x$ with respect to the inner unit normal, arranged in non-decreasing order: $\eta_1(x)\leq\ldots\leq\eta_n(x)$. For each $p$ ($1\leq p \leq n$), let the $p$-curvature $\sigma_p(x)$ be the sum of the $p$ smallest principal curvatures, i.e., $\sigma_p(x)=\eta_{1}(x)+\ldots+\eta_{p}(x)$. We then define $$\sigma_p(\partial M)=\inf_{x\in\partial M}\sigma_p(x).$$ 
In particular, from \eqref{2-Sec-5} we have immediately
\begin{equation}\label{2-thm-formula-1}
\langle S^{[p]}\omega,\omega\rangle\geq\sigma_p(\partial M)|\omega|^2.
\end{equation}

In order to state our results, we let $\lambda'_{1,p}(\partial M)$ (resp. $\lambda''_{1,p}(\partial M)$) be the first eigenvalue of the Hodge Laplacian when restricted to exact (resp. co-exact) $p$-forms of $\partial M$, and let $\lambda_{1,p}(\partial M)$ be the first positive eigenvalue of $\Delta_f^H$ on $\partial M$, respectively. Then the following Hodge decomposition theorem
$$A^p=H^p\oplus dA^{p-1}\oplus\delta_f A^{p+1}$$
implies
$$\lambda_{1,p}(\partial M)=\min\{\lambda'_{1,p}(\partial M),\lambda''_{1,p}(\partial M)\}$$
and
\begin{equation}\label{dual-1}
\begin{cases}
\lambda''_{1,p}(\partial M)=\lambda'_{1,p+1}(\partial M),\\
\lambda''_{1,p}(\partial M)=\lambda'_{1,n-p}(\partial M).
\end{cases}
\end{equation}
For the proof of \eqref{dual-1}, see Lemma \ref{1-App-Lema-3} in Appendix. Taking $V=1$ in \eqref{1-Th-Formula-1}, we obtain the following

\begin{thm}\label{1-Th-2}
Let $(M, g, d\mu)$ be a compact smooth metric measure space with $W^{[p]}_{f}\geq 0$ and $\sigma_p(\partial M)>0$ on $\partial M$. If $\sigma_{n-p+1}(\partial M)+\inf(f_N)>0$, then we have
\begin{equation}\label{2-Th2-1}
\lambda'_{1,p}(\partial M)\geq \sigma_p(\partial M)[\sigma_{n-p+1}(\partial M)+\inf(f_N)].
\end{equation}
\end{thm}

We remark that if $f$ is a constant, then \eqref{2-Th2-1} becomes $\lambda'_{1,p}(\partial M)\geq \sigma_p(\partial M)\sigma_{n-p+1}(\partial M)$. In this case, it was shown in \cite{Savo2011} that if $M$ has non-negative Ricci curvature, then the equality holds if and only if $M$ is a Euclidean ball. However, for non-trivial function $f$, we can not give a characterization in \eqref{2-Th2-1} for equality case.

Recall the following well-known facts: The absolute cohomology space $H^p(M,\mathbb{R})$ of $M$ in degree $p$
with real coefficients, is isomorphic to the space of harmonic $p$-forms satisfying the boundary problem
\begin{equation}\label{3-Th-Formula-1}
\begin{cases}
d\widehat{\omega}=\delta_f\widehat{\omega}=0,\ \  {\rm on}\ M;\\
i_N\widehat{\omega}=0,\ \ \ \ \ \ \ \ \ {\rm on}\ \partial M.
\end{cases}
\end{equation}
Then, we have

\begin{thm}\label{1-Th-3}
Let $(M, g, d\mu)$ be a compact smooth metric measure space with $W^{[p]}_{f,V}\geq 0$ and $V$ be a positive smooth function defined on $M$.

$(i)$ If $\sigma_p(\partial M)>(\ln V)_N$, then $H^p(M,\mathbb{R})=0$.

$(ii)$ If $\sigma_p(\partial M)\geq(\ln V)_N$ and $H^p(M,\mathbb{R})\neq0$, then $M$ admits a non-trivial parallel $p$-form and $\sigma_p(\partial M)=(\ln V)_N$ for all $x\in\partial M$.
\end{thm}

We remark that if $V$ is a constant, then our Theorem \ref{1-Th-3} becomes Theorem 4 in \cite{Savo2011}.

Next, we consider the Steklov eigenvalue problem for differential forms. Given an integer $0\leq p\leq n$ and $\omega\in A^p(\partial M)$, there exists a solution $\widehat{\omega}\in A^p(M)$ satisfying
\begin{equation}\label{5-Th-formula-1}
\begin{cases}
\Delta_f^H\widehat{\omega}=0,\ \delta_f\widehat{\omega}=0,\ {\rm on}\ M;\\
J^*\widehat{\omega}=\omega,\ \ \ \ \ \ \ \ \ \  \ \ \ \  \ {\rm on}\ \partial M.
\end{cases}
\end{equation}
The Dirichlet-to-Neumann operator $\mathcal{D}^{[p]}$ defined by
$$\mathcal{D}^{[p]}(\omega):=-i_Nd\widehat{\omega}$$
is to study the spectral properties of $\mathcal{D}^{[p]}$, which can be restricted to the subspace of co-closed $p$-forms in $A^p(\partial M)$, see \cite{Karpukhin2019}. Specially,
\begin{equation}\label{5-Th-formula-2}
-i_Nd\widehat{\omega}=\sigma \omega.
\end{equation}
In this case, the restricted operator admits a sequence of eigenvalues
$$0\leq\sigma_1^{[p]}(\partial M)\leq\sigma_2^{[p]}(\partial M)\leq\ldots\nearrow\infty.$$

\begin{thm}\label{1-Th-5}
Let $(M, g, d\mu)$ be a compact smooth metric measure space with the non-negative $(p+1)$-st $f$-Weitzenb\"{o}ck curvature and sectional curvature, where $1\leq p\leq n$.  Denote by $\sigma(\partial M)$ the first non-zero eigenvalue of the Dirichlet-to-Neumann operator $\mathcal{D}^{[p]}$. If the principal curvatures of $\partial M$ are no less than $c>0$, then
\begin{equation}\label{1-Th-5-formula-1}
\sigma(\partial M)\geq (p+1)c.
\end{equation}
\end{thm}

When $f$ is constant, it has been proved in \cite{Xiong2024} that the equality holds for a Euclidean ball with the radius $\frac{1}{c}$. In particular, for $p=0$, the result in Theorem \ref{1-Th-5} is related to the well-known Escobar's conjecture
(for more details, we refer the reader to \cite{Escobar1997,Escobar1999,XX2024} and the references therein).

Note that
$$
\sigma^{[p]}_k(\partial M)=\inf_{E_k}\sup_{\widehat{\omega}\in E_k\backslash\{0\}}\frac{\int_M|d\widehat{\omega}|^2\,e^{-f}dv}{\int_{\partial M}|\omega|^2\,e^{-f}da},
$$
where $E_k$ runs over all $k$-dimensional subspaces in $A^p(M)$ and $J^*E_k$  is included in co-closed forms, see the formula (2.6) in \cite{Xiong2024}. Then, we have

\begin{thm}\label{1-Th-6}
Assumptions are the same as in Theorem \ref{1-Th-5} except we suppose $1\leq p\leq n-1$. Denote by $\sigma^{[p]}_k(\partial M)$ the eigenvalues of the operator $\mathcal{D}^{[p]}$ and $\lambda^{[p]}_k(\partial M)$ the eigenvalues of the weighted Hodge Laplacian with respect to the induced metric on $\partial M$, both operators being restricted to the co-closed differential forms, where $k\geq1$. If $(n-p)c-\sup|\nabla f|>0$, then we have
\begin{equation}\label{1-Th-6-formula-1}
\sigma^{[p]}_k(\partial M)\leq\frac{1}{(n-p)c-\sup|\nabla f|}\lambda^{[p]}_k(\partial M).
\end{equation}
\end{thm}

When $f$ is constant, \eqref{1-Th-6-formula-1} becomes $\sigma^{[p]}_k(\partial M)\leq\frac{1}{(n-p)c}\lambda^{[p]}_k$, which has been proved by Xiong in \cite{Xiong2024}. For some other interesting results for comparisons between the Steklov eigenvalues and the Laplacian eigenvalues on the boundary, see \cite{Colbois2020,HMZ2023,Karpukhin2017,Kwong2016,Xiong2018,Wang2009,Wang2018} and the references therein.

The spectrum of $\Delta_f^H$ acting on $p$-forms consists of a non-decreasing, unbounded sequence of eigenvalues
with finite multiplicities
$${\rm Spec}(\Delta_f^H)=\{0\leq\lambda_1^{(p)}\leq\lambda_2^{(p)}\leq\ldots\lambda_j^{(p)}\leq\ldots\}.$$

Next, we consider closed submanifolds in a weighted Euclidean space. By applying the weighted Bochner formula \eqref{2-Sec-4},
we derive universal inequalities for the sum of eigenvalues with respect to the weighted Hodge Laplacian, which can be seen as a generalization to Levitin-Parnovski inequality (see \cite{LP2002}).

\begin{thm}\label{1-Th-7}
Let $X:(M^m,g)\stackrel{{\rm can}}{\rightarrow}(\mathbb{R}^{M}, e^{-f}\delta)$ be a closed isometric immersion with the mean curvature vector field $\mathbf{H}$, where $\delta$ denotes the Euclidean metric. Then, for any $p\in\{1,\ldots,m\}$ and $j\in\mathbb{N}$, we have
\begin{align}\label{7-thm-1}
\sum_{l=1}^m\lambda_{j+l}^{(p)}\leq&(m+4)\lambda_j^{(p)}+2\sup\Delta_ff\notag\\
&+\int_M[(m^2|\mathbf{H}|^2+|\nabla f|^2)|\omega_j|^2-4\langle W^{[p]}_f(\omega_j),\omega_j\rangle]\,e^{-f}dv,
\end{align}
where $\{\lambda_j\}_{j=1}^\infty$ are the eigenvalues of $\Delta_f^H$ and $\{\omega_j\}_{j=1}^\infty$ is a corresponding orthonormal basis of $p$-eigenforms.
\end{thm}

For some generalizations of the Levitin-Parnovski inequality with respect to Laplacian, see \cite{CC2008,CHW2010,Sun2008} and the references therein. When $f$ is constant, our Theorem \ref{1-Th-7} becomes Theorem 2.1 of Ilias and Makhoul in \cite{Ilias2012}. On the other hand,
$$m^2|\mathbf{H}|^2+|\nabla f|^2=|m\mathbf{H}-(\overline{\nabla} f)^{\top}|^2=|m\mathbf{H}_f-\overline{\nabla} f|^2,$$
where $\mathbf{H}_f=\mathbf{H}+\frac{1}{m}(\overline{\nabla} f)^{\bot}$,  $\overline{\nabla}$ is the gradient operator with respect to the Euclidean metric $\delta$, $\top$ and $\bot$ are the projection onto the tangent bundle and normal bundle of $X$, respectively. Hence, \eqref{7-thm-1} can be written as
\begin{align}\label{7-thm-2}
\sum_{l=1}^m\lambda_{j+l}^{(p)}\leq&(m+4)\lambda_j^{(p)}+2\sup\Delta_ff\notag\\
&+\int_M[|m\mathbf{H}_f-\overline{\nabla} f|^2|\omega_j|^2-4\langle W^{[p]}_f(\omega_j),\omega_j\rangle]\,e^{-f}dv.
\end{align}
In particular, under some suitable conditions, some related conclusions as in \cite{Ilias2012} can also be obtained.

This paper is organized as follows: In Section 2, we provide the necessary weighted Bochner formula and related Reilly type formula with respect to weighted Hodge Laplacian; In Section 3, we study compact manifolds with boundary and give proofs of Theorems \ref{1-Th-2}-\ref{1-Th-6}, respectively; In Section 4, Levitin-Parnovski eigenvalue inequality on closed submanifolds of a weighted Euclidean space has been proved; In Section 5, some necessary formulas are provided.

\section{The weighted Bochner formula and related Reilly type formula}
In this section, we present two fundamental propositions that are crucial for our subsequent proofs.
In order to state our results, we adopt the following conceptions in
\cite{IRS2020,Xiong2024,Savo2011}:
For any $\omega\in A^p(M)$, we have
\begin{equation}\label{1-Boch-1}
(i_N\omega)(X_1,\ldots,X_{p})=\omega(N,X_1,\ldots,X_{p});
\end{equation}
\begin{equation}\label{1-Boch-2}
(\delta\omega)(X_1,\ldots,X_{p-1})=-(\nabla_{e_A}\omega)(e_A,X_1,\ldots,X_{p-1});
\end{equation}
\begin{equation}\label{1-Boch-3}
d\omega=e_A^*\wedge \nabla_{e_A}\omega,
\end{equation}
where $\{e_A\}_{A=1}^{n+1}$ is a local orthonormal frame of $M$, $e_A^*$ is the dual 1-form of $e_A$ and
$X_r\in\mathcal{X}(M)$, $r=1,2,\ldots,p$. In particular, it is easy to see from \eqref{1-Boch-1} that \eqref{1-Boch-2} can be written as
\begin{equation}\label{1-Boch-222}
\delta\omega=-i_{e_A}(\nabla_{e_A}\omega).
\end{equation}

The following proposition, which is known in the literature (see \cite[Proposition 4.4]{Petersen2012}), is stated here using our notation for completeness and to set the stage for our main results.
\begin{prop}\cite{Petersen2012}\label{1-sec-Prop-1}
Let $(M, g, d\mu)$ be a smooth metric measure space. For the weighted Hodge Laplacian given by \eqref{2-Sec-3}, we have
the following Bochner formula:
\begin{equation}\label{1-lem-1}
\Delta_f^H=\nabla_f^*\nabla+W_f^{[p]},
\end{equation}
where
\begin{equation}\label{1-lem-111}
\nabla_f^*\nabla=\nabla^*\nabla+\nabla_{\nabla f},\ \ \ \
 W_f^{[p]}=W^{[p]}+(\nabla^2f)^{[p]}.
 \end{equation}
In particular, for any $\omega\in A^p(M)$, we have
\begin{equation}\label{1-lem-222}
\frac{1}{2}\Delta_f|\omega|^2=\langle\Delta_f^H\omega,\omega\rangle-\langle W^{[p]}_f(\omega),\omega\rangle-|\nabla\omega|^2.
\end{equation}

\end{prop}

\proof
From the definition of $\delta_f$ in \eqref{2-Sec-2}, we have
$$\aligned
(\delta_f\omega)(X_1,\ldots,X_{p-1})=&e^{f}\delta(e^{-f}\omega)(X_1,\ldots,X_{p-1})\\
=&-e^{f}\nabla_{e_A}(e^{-f}\omega)(e_A,X_1,\ldots,X_{p-1})\\
=&-e^{f}(-e^{-f}f_A\omega+e^{-f}\nabla_{e_A}\omega)(e_A,X_1,\ldots,X_{p-1})\\
=&(f_A\omega-\nabla_{e_A}\omega)(e_A, X_1,\ldots,X_{p-1})\\
=&(i_{\nabla f}\omega)(X_1,\ldots,X_{p-1})+(\delta\omega)(X_1,\ldots,X_{p-1})
\endaligned$$
which yields
\begin{equation}\label{1-lem-2}
\delta_f=\delta+i_{\nabla f}.
\end{equation}
Therefore, it follows from \eqref{1-lem-2} that
\begin{align}\label{1-lem-3}
\Delta_f^H&=d\delta_f+\delta_fd\notag\\
&=d(\delta+i_{\nabla f})+(\delta+i_{\nabla f})d\notag\\
&=\Delta^H+di_{\nabla f}+i_{\nabla f}d.
\end{align}
For any fixed index $\alpha$,
\begin{align*}
&[\nabla_{X_\alpha}(i_{\nabla f}\omega)](X_1,\ldots,\widehat{X_{\alpha}},\ldots,X_p)\notag\\
=&\nabla_{X_{\alpha}}[(i_{\nabla f}\omega)(X_1,\ldots,\widehat{X_{\alpha}},\ldots,X_p)]\notag\\
&-(i_{\nabla f}\omega)(\nabla_{X_{\alpha}}X_1,\ldots,\widehat{X_{\alpha}},\ldots,X_p)-\ldots-(i_{\nabla f}\omega)(X_1,\ldots,\widehat{X_{\alpha}},\ldots,\nabla_{X_{\alpha}}X_p)\notag\\
=&X_{\alpha}[\omega(\nabla f,X_1,\ldots,\widehat{X_{\alpha}},\ldots,X_p)]\notag\\
&-\omega(\nabla f,\nabla_{X_{\alpha}}X_1,\ldots,\widehat{X_{\alpha}},\ldots,X_p)-\ldots-\omega(\nabla f,X_1,\ldots,\widehat{X_{\alpha}},\ldots,\nabla_{X_{\alpha}}X_p)\notag\\
=&(\nabla_{X_{\alpha}}\omega)(\nabla f,X_1,\ldots,\widehat{X_{\alpha}},\ldots,X_p)+\omega(\nabla_{X_{\alpha}}\nabla f,X_1,\ldots,\widehat{X_{\alpha}},\ldots,X_p)\notag\\
=&(\nabla\omega)(X_{\alpha},\nabla f,X_1,\ldots,\widehat{X_{\alpha}},\ldots,X_p)+\omega(\nabla_{X_{\alpha}}\nabla f,X_1,\ldots,\widehat{X_{\alpha}},\ldots,X_p),
\end{align*}
where $\widehat{X_{\alpha}}$ means that the argument $X_{\alpha}$ disappears, this shows
\begin{align}\label{1-lem-4}
\sum_{\alpha=1}^p(-1)^{\alpha+1}&[\nabla_{X_{\alpha}}(i_{\nabla f}\omega)](X_1,\ldots,\widehat{X_{\alpha}},\ldots,X_p)\notag\\
&+\sum_{\alpha=1}^p(-1)^\alpha(\nabla\omega)(X_{\alpha},\nabla f,X_1,\ldots,\widehat{X_{\alpha}},\ldots,X_p)\notag\\
=&\sum_{\alpha=1}^p\omega(X_1,\ldots,\nabla_{X_{\alpha}}\nabla f,\ldots,X_p)\notag\\
=&\sum_{\alpha=1}^p\omega(X_1,\ldots,\nabla^2 f(X_{\alpha}),\ldots,X_p)\notag\\
=&[(\nabla^2 f)^{[p]}\omega](X_1,\ldots,X_p).
\end{align}
From \eqref{1-Boch-3}, it holds that
\begin{align*}\label{1-lem-55}
(&di_{\nabla f}\omega+i_{\nabla f}d\omega)(X_1,\ldots,X_p)\notag\\
=&[d(i_{\nabla f}\omega)](X_1,\ldots,X_p)+d\omega(\nabla f,X_1,\ldots,X_p)\notag\\
=&e_A^*\wedge[\nabla_{e_A}(i_{\nabla f}\omega)](X_1,\ldots,X_p)+e_A^*\wedge(\nabla_{e_A}\omega)(\nabla f,X_1,\ldots,X_p)\notag\\
=&\sum_{\alpha=1}^p(-1)^{\alpha+1}[\nabla_{X_\alpha}(i_{\nabla f}\omega)](X_1,\ldots,\widehat{X_\alpha},\ldots,X_p)\notag\\
&+(\nabla\omega)(\nabla f,X_1,\ldots,X_p)+\sum_{\alpha=1}^p(-1)^\alpha(\nabla\omega)(X_\alpha,\nabla f,X_1,\ldots,\widehat{X_\alpha},\ldots,X_p)\notag\\
=&(\nabla_{\nabla f}\omega)(X_1,\ldots,X_p)+[(\nabla^2 f)^{[p]}\omega](X_1,\ldots,X_p),
\end{align*}
where in the last equality, we used \eqref{1-lem-4}, which implies
\begin{equation}\label{1-lem-5}
di_{\nabla f}+i_{\nabla f}d=\nabla_{\nabla f}+(\nabla^2 f)^{[p]}.
\end{equation}
Therefore, \eqref{1-lem-3} can be written as
$$\aligned
\Delta_f^H=&\Delta^H+\nabla_{\nabla f}+(\nabla^2 f)^{[p]}\\
=&\nabla_f^*\nabla+W_f^{[p]},
\endaligned$$
where we used the following well-known Bochner formula (see \cite{Savo2011,Peter1996}):
\begin{align}\label{1-lem-6}
\Delta^H=\nabla^*\nabla+W^{[p]}.
\end{align}

On the other hand, by the definition of weighted Laplacian, we have
$$\aligned
\frac{1}{2}\Delta_f|\omega|^2=&-\frac{1}{2}e^f(e^{-f}|\omega|^2_{A})_{A}\\
=&-e^{f}(e^{-f}\omega_I\omega_{I,A})_{,A}\\
=&f_A\omega_{I,A}\omega_I-|\omega _{I,A}|^2-\omega_I\omega_{I,AA}\\
=&\langle\nabla_{\nabla f}\omega,\omega\rangle-|\nabla\omega|^2+\langle\nabla^*\nabla\omega,\omega\rangle\\
=&\langle\nabla^*_f\nabla\omega,\omega\rangle-|\nabla\omega|^2\\
=&\langle\Delta_f^H\omega,\omega\rangle-\langle W^{[p]}_f(\omega),\omega\rangle-|\nabla\omega|^2,
\endaligned$$
where $|\omega|^2=\sum\omega _{I}^2$ with multi-index $I=(i_1,i_2,\ldots,i_p)$, and in the last equality we used \eqref{1-lem-1}. Therefore, we complete the proof of Proposition \ref{1-sec-Prop-1}.

With the help of formula \eqref{1-lem-222} in Proposition \ref{1-sec-Prop-1}, we give the following weighted Reilly type formula associated with the weighted measure:

\begin{prop}\label{prop-2}
Let $(M, g, d\mu)$ be a compact smooth metric measure space and $V$ be a smooth function defined on $M$. Then,
\begin{align}\label{2-prop-1}
\int_M&V(|\delta_f\omega|^2+|d\omega|^2-|\nabla\omega|^2)\,e^{-f}dv\notag\\
=&\int_M[-2\langle\omega,i_{\nabla V}(d\omega)\rangle +V\langle W_{f,V}^{[p]}(\omega),\omega\rangle]\,e^{-f}dv\notag\\
&+\int_{\partial M}[-V_N|J^*\omega|^2+2V\langle\delta^{\partial M}_f(J^*\omega),i_N\omega\rangle+ V\mathcal{B}_f(\omega,\omega)]\,e^{-f}da,
\end{align}
where $W_{f,V}^{[p]}$ is given by \eqref{1-Th-Formula-2} and
\begin{equation}\label{2-prop-2}
\mathcal{B}_f(\omega,\omega)=\langle S^{[p]}(J^*\omega),J^*\omega\rangle+nH_f|i_N\omega|^2-\langle S^{[p-1]}(i_N\omega),i_N\omega\rangle,
\end{equation}
which is  equivalent to
\begin{equation}\label{2-prop-3}
\mathcal{B}_f(\omega,\omega)=\langle S^{[p]}(J^*\omega),J^*\omega\rangle+\langle S^{[n+1-p]}(J^* *\omega),J^* *\omega\rangle+f_N|i_N\omega|^2.
\end{equation}
\end{prop}

\proof
Multiplying both sides of \eqref{1-lem-222} by $V$ and integrating it, we have
\begin{align}\label{1-prop-3}
\frac{1}{2}\int_MV\Delta_f|\omega|^2=&\int_MV(\langle\Delta_f^H\omega,\omega\rangle-\langle W^{[p]}_f(\omega),\omega\rangle-|\nabla\omega|^2)\,e^{-f}dv.
\end{align}
Since
$$\aligned
\Delta_f(V|\omega|^2)=&|\omega|^2\Delta_fV+V\Delta_f|\omega|^2-2\langle \nabla V,\nabla|\omega|^2\rangle\\
=&|\omega|^2\Delta_fV+V\Delta_f|\omega|^2-4\langle\omega,\nabla_{\nabla V}\omega\rangle,
\endaligned$$
then \eqref{1-prop-3} can be written as
\begin{align}\label{1-prop-4}
\int_M&\Big[\frac{1}{2}|\omega|^2\Delta_fV+V(\langle\Delta_f^H\omega,\omega\rangle-\langle W^{[p]}_f(\omega),\omega\rangle-|\nabla\omega|^2)-2\langle\omega,\nabla_{\nabla V}\omega\rangle\Big]\,e^{-f}dv\notag\\
=&\frac{1}{2}\int_M\Delta_f(V|\omega|^2)\,e^{-f}dv\notag\\
=&\frac{1}{2}\int_{\partial M}(V_N|\omega|^2+V|\omega|^2_N)\,e^{-f}da.
\end{align}
Note that
\begin{align}\label{1-prop-4444}
\delta_f(V\omega)=&(\delta+i_{\nabla f})(V\omega)\notag\\
=&-i_{e_A}\nabla_{e_A}(V\omega)+Vi_{\nabla f}\omega\notag\\
=&-i_{e_A}(V_A\omega+V\nabla_{e_A}\omega)+Vi_{\nabla f}\omega\notag\\
=&-i_{\nabla V}\omega+V\delta_f\omega.
\end{align}
Thus, applying \eqref{1-App-Lema-1}, we have
\begin{align}\label{1-prop-9}
\int_MV\langle d\delta_f\omega,\omega\rangle\,e^{-f}dv=&\int_M\langle \delta_f\omega,\delta_f(V\omega)\rangle\,e^{-f}dv
-\int_{\partial M}\langle J^*(\delta_f\omega),i_N(V\omega)\rangle\,e^{-f}da\notag\\
=&\int_M\langle \delta_f\omega,-i_{\nabla V}\omega+V\delta_f\omega\rangle\,e^{-f}dv
-\int_{\partial M}V\langle J^*(\delta_f\omega),i_N\omega\rangle\,e^{-f}da\notag\\
=&\int_M(-\langle \delta_f\omega,i_{\nabla V}\omega\rangle+V|\delta_f\omega|^2)\,e^{-f}dv
-\int_{\partial M}V\langle J^*(\delta_f\omega),i_N\omega\rangle\,e^{-f}da\notag\\
=&\int_M(-\langle d(i_{\nabla V}\omega),\omega\rangle+V|\delta_f\omega|^2)\,e^{-f}dv-\int_{\partial M}(\langle J^*(i_{\nabla V}\omega),i_N\omega\rangle\notag\\
&+V\langle J^*(\delta_f\omega),i_N\omega\rangle)\,e^{-f}da.
\end{align}
For a Lipschitz vector field $F$, it has been proved in \cite{Xiong2024} that
\begin{equation}\label{1-prop-5}
d(i_F\omega)=-i_F(d\omega)+\nabla_F\omega+\nabla F(\omega).
\end{equation}
Thus, we have
$$d(i_{\nabla V}\omega)=-i_{\nabla V}(d\omega)+\nabla_{\nabla V}\omega+(\nabla^2 V)^{[p]}\omega,$$
and \eqref{1-prop-9} becomes
\begin{align}\label{1-prop-10}
\int_MV\langle d\delta_f\omega,\omega\rangle e^{-f}dv=&\int_M[\langle i_{\nabla V}(d\omega)-\nabla_{\nabla V}\omega-(\nabla^2 V)^{[p]}\omega,\omega\rangle+V|\delta_f\omega|^2]\,e^{-f}dv\notag\\
&-\int_{\partial M}(\langle J^*(i_{\nabla V}\omega),i_N\omega\rangle+V\langle J^*(\delta_f\omega),i_N\omega\rangle) \,e^{-f}da.
\end{align}

On the other hand, applying \eqref{1-App-Lema-1} again yields
\begin{align}\label{1-prop-11}
\int_MV\langle \delta_fd\omega,\omega\rangle e^{-f}dv=&\int_M\langle d\omega,d(V\omega)\rangle e^{-f}dv+\int_{\partial M}\langle J^*(V\omega),i_Nd\omega\rangle e^{-f}da\notag\\
=&\int_M[\langle d\omega,dV\wedge\omega\rangle+\langle d\omega,Vd\omega\rangle]\,e^{-f}dv+\int_{\partial M}V\langle J^*\omega,i_Nd\omega\rangle e^{-f}da\notag\\
=&\int_M[\langle\omega,\delta_f(dV\wedge\omega)\rangle+V|d\omega|^2]\,e^{-f}dv+\int_{\partial M}[-\langle J^*\omega,i_N(dV\wedge\omega)\rangle\notag\\
&+V\langle J^*\omega,i_Nd\omega\rangle]\,e^{-f}da.
\end{align}
Noticing
\begin{align*}\label{1-prop-12}
&i_{\nabla f}(dV\wedge\omega)(X_1,\ldots,X_p)\notag\\
=&(dV\wedge\omega)(\nabla f,X_1,\ldots,X_p)\notag\\
=&\langle \nabla V,\nabla f\rangle\omega(X_1,\ldots,X_p)+\sum_{\alpha=1}^p(-1)^\alpha\langle\nabla V,X_\alpha\rangle(i_{\nabla f}\omega)(X_1,\ldots,\widehat{X_\alpha},\ldots X_p)\notag\\
=&(\langle \nabla V,\nabla f\rangle\omega-dV\wedge i_{\nabla f}\omega)(X_1,\ldots,X_p),
\end{align*}
we have
\begin{equation}\label{1-prop-13}
i_{\nabla f}(dV\wedge\omega)=\langle \nabla V,\nabla f\rangle \omega-dV\wedge i_{\nabla f}\omega.
\end{equation}
Thus, with the help of the following(see Proposition 25 in \cite{Xiong2024}) formula:
\begin{equation}\label{1-prop-1000}
\delta(dV\wedge\omega)=(\Delta V)\omega-\nabla_{\nabla V}\omega+(\nabla^2 V)^{[p]}\omega-dV\wedge\delta\omega,
\end{equation}
we obtain
\begin{align}\label{1-prop-2000}
\delta_f(dV\wedge\omega)=&(\delta+i_{\nabla f})(dV\wedge\omega)\notag\\
=&(\Delta_f V)\omega-\nabla_{\nabla V}\omega+(\nabla^2 V)^{[p]}\omega-dV\wedge\delta_f\omega
\end{align}
and \eqref{1-prop-11} becomes
\begin{align}\label{1-prop-14}
\int_M&V\langle \delta_fd\omega,\omega\rangle e^{-f}dv\notag\\
=&\int_M[(\Delta_f V)|\omega|^2+\langle\omega,-\nabla_{\nabla V}\omega+(\nabla^2 V)^{[p]}\omega\rangle+V|d\omega|^2\notag\\
&-V|\delta_f\omega|^2+V\langle d\delta_f\omega,\omega \rangle]\,e^{-f}dv+\int_{\partial M}[V\langle J^*(\delta_f\omega),i_N\omega\rangle \notag\\
&-\langle J^*\omega,i_N(dV\wedge\omega)\rangle+ V\langle J^*\omega,i_Nd\omega\rangle]\,e^{-f}da,
\end{align}
where we used
\begin{align}\label{1-prop-16}
&-\int_M\langle\omega,dV\wedge\delta_f\omega\rangle \,e^{-f}dv\notag\\
=&-\int_M\langle\omega,d(V\delta_f\omega)-Vd\delta_f\omega\rangle\, e^{-f}dv\notag\\
=&\int_M[-V|\delta_f\omega|^2+V\langle d\delta_f\omega,\omega \rangle]\,e^{-f}dv+\int_{\partial M} V\langle J^*(\delta_f\omega),i_N\omega\rangle e^{-f}da.
\end{align}
As a result, by using \eqref{1-prop-10} and \eqref{1-prop-14}, we have
\begin{align}\label{1-prop-17}
\int_M&V\langle\Delta_f^H\omega,\omega\rangle\, e^{-f}dv\notag\\
=&\int_MV\langle d\delta_f\omega,\omega\rangle e^{-f}dv+\int_MV\langle \delta_fd\omega,\omega\rangle\,e^{-f}dv\notag\\
=&\int_M[2V\langle d\delta_f\omega,\omega\rangle+(\Delta_f V)|\omega|^2+\langle\omega,-\nabla_{\nabla V}\omega+(\nabla^2V)^{[p]}\omega\rangle\notag\\
&+V|d\omega|^2-V|\delta_f\omega|^2]\,e^{-f}dv+\int_{\partial M}[V\langle J^*(\delta_f\omega),i_N\omega\rangle \notag\\
&-\langle J^*\omega,i_N(dV\wedge\omega)\rangle+V\langle J^*\omega,i_Nd\omega\rangle]\, e^{-f}da\notag\\
=&\int_M[2\langle i_{\nabla V}(d\omega),\omega\rangle-3\langle\nabla_{\nabla V}\omega,\omega\rangle-\langle(\nabla^2V)^{[p]}\omega,\omega\rangle+V|\delta_f\omega|^2\notag\\
&+|\omega|^2\Delta_fV+V|d\omega|^2]\,e^{-f}dv+\int_{\partial M}[-2\langle J^*(i_{\nabla V}\omega),i_N\omega\rangle\notag\\
&-V\langle J^*(\delta_f\omega),i_N\omega\rangle-\langle J^*\omega,i_N(dV\wedge\omega)\rangle+V\langle J^*\omega,i_Nd\omega\rangle]\,e^{-f}da,
\end{align}
and substituting \eqref{1-prop-17} into \eqref{1-prop-4}, we obtain
\begin{align}\label{1-prop-18}
\frac{1}{2}&\int_{\partial M}(V_N|\omega|^2+V|\omega|^2_N)\,e^{-f}da\notag\\
=&\int_M[2\langle i_{\nabla V}(d\omega),\omega\rangle-5\langle\nabla_{\nabla V}\omega,\omega\rangle-\langle(\nabla^2V)^{[p]}\omega,\omega\rangle\notag\\
&+V|\delta_f\omega|^2+\frac{3}{2}|\omega|^2\Delta_fV-V\langle W^{[p]}_f(\omega),\omega\rangle-V|\nabla\omega|^2+V|d\omega|^2]\,e^{-f}dv\notag\\
&+\int_{\partial M}[-2\langle J^*(i_{\nabla V}\omega),i_N\omega\rangle-V\langle J^*(\delta_f\omega),i_N\omega\rangle-\langle J^*\omega,i_N(dV\wedge\omega)\rangle\notag\\
&+V\langle J^*\omega,i_Nd\omega\rangle]e^{-f}da.
\end{align}
Applying the divergence theorem, we have
\begin{align}\label{1-prop-19}
\int_M\langle\nabla_{\nabla V}\omega,\omega\rangle \,e^{-f}dv=&\frac{1}{2}\int_M\langle\omega,\omega\rangle_{\nabla V}\,e^{-f}dv\notag\\
=&-\frac{1}{2}\int_{\partial M}|\omega|^2V_Ne^{-f}da+\frac{1}{2}\int_{M}|\omega|^2\Delta_fV\,e^{-f}dv
\end{align}
and \eqref{1-prop-18} becomes
\begin{align}\label{1-prop-20}
&\int_{\partial M}[\frac{1}{2}V|\omega|^2_N-2V_N|\omega|^2]\,e^{-f}da\notag\\
=&\int_M[2\langle i_{\nabla V}(d\omega),\omega\rangle-|\omega|^2\Delta_fV-\langle(\nabla^2V)^{[p]}\omega,\omega\rangle+V|\delta_f\omega|^2\notag\\
&-V\langle W^{[p]}_f(\omega),\omega\rangle-V|\nabla\omega|^2+V|d\omega|^2]\,e^{-f}dv\notag\\
&+\int_{\partial M}[-2\langle J^*(i_{\nabla V}\omega),i_N\omega\rangle-V\langle J^*(\delta_f\omega),i_N\omega\rangle-\langle J^*\omega,i_N(dV\wedge\omega)\rangle\notag\\
&+V\langle J^*\omega,i_Nd\omega\rangle]\,e^{-f}da,
\end{align}
which implies
\begin{align}\label{1-prop-21}
&\int_{\partial M}V[\langle J^*(\nabla_N\omega),J^*\omega\rangle
+\langle i_N(\nabla_N\omega),i_N\omega\rangle-2V_N|\omega|^2]\,e^{-f}da\notag\\
=&\int_M[2\langle i_{\nabla V}(d\omega),\omega\rangle-|\omega|^2\Delta_fV-\langle\nabla^2V(\omega),\omega\rangle+V|\delta_f\omega|^2-V\langle W^{[p]}_f(\omega),\omega\rangle\notag\\
&-V|\nabla\omega|^2+V|d\omega|^2]\,e^{-f}dv\notag\\
&+\int_{\partial M}[-2\langle J^*(i_{\nabla V}\omega),i_N\omega\rangle-V\langle J^*(\delta_f\omega),i_N\omega\rangle-\langle J^*\omega,V_NJ^*\omega-d^{\partial M}V\wedge i_N\omega\rangle\notag\\
&+V\langle J^*\omega,i_Nd\omega\rangle]\,e^{-f}da,
\end{align}
where we used
$$\frac{1}{2}|\omega|^2_N=\langle\nabla_N\omega,\omega\rangle=\langle J^*(\nabla_N\omega),J^*\omega\rangle
+\langle i_N(\nabla_N\omega),i_N\omega\rangle$$
and
$$i_N(dV\wedge\omega)=V_NJ^*\omega-d^{\partial M}V\wedge i_N\omega.$$
According to \eqref{3-Appx2-2} and \eqref{3-Appx2-3}, then \eqref{1-prop-21} becomes
\begin{align}\label{1-prop-23}
0=&\int_M[2\langle i_{\nabla V}(d\omega),\omega\rangle-V\langle W_{f,V}^{[p]}(\omega),\omega\rangle
+V(|\delta_f\omega|^2-|\nabla\omega|^2+|d\omega|^2)]\,e^{-f}dv\notag\\
&+\int_{\partial M}[-V\langle d^{\partial M}(i_N\omega),J^*\omega\rangle -2\langle J^*(i_{\nabla V}\omega),i_N\omega\rangle \notag\\
&+2|\omega|^2V_N-V\mathcal{B}_f(\omega,\omega)
-V\langle\delta^{\partial M}_f(J^*\omega),i_N\omega\rangle \notag\\
&-|J^*\omega|^2V_N+\langle J^*\omega,d^{\partial M}V\wedge i_N\omega\rangle]\,e^{-f}da,
\end{align}
where $\mathcal{B}_f$ is given by \eqref{1-Th-Formula-3}.

On the boundary $\partial M$, we have
$$|\omega|^2=|J^*\omega|^2+|i_N\omega|^2,\ \ \ i_{\nabla V}\omega=i_{\nabla^{\partial M} V}\omega+V_Ni_N\omega$$
and then
\begin{align}\label{1-prop-24}
&\int_{\partial M}[-2\langle J^*(i_{\nabla V}\omega),i_N\omega\rangle +2|\omega|^2V_N-|J^*\omega|^2V_N]\,e^{-f}da\notag\\
&=\int_{\partial M}[-2\langle J^*(i_{\nabla^{\partial M}V}\omega),i_N\omega\rangle+|J^*\omega|^2V_N]\,e^{-f}da.
\end{align}
In addition, by applying \eqref{3-Appx-1} and \eqref{1-prop-4444}, we have
\begin{align}\label{1-prop-25}
-&\int_{\partial M}V\langle d^{\partial M}i_N\omega,J^*\omega\rangle\,e^{-f}da
=-\int_{\partial M}\langle i_N\omega,\delta_f^{\partial M}(VJ^*\omega)\rangle \,e^{-f}da\notag\\
=&\int_{\partial M}[\langle i_N\omega,i_{\nabla^{\partial M}V}(J^*\omega)\rangle-V\langle i_N\omega, \delta_f^{\partial M}(J^*\omega)\rangle]\,e^{-f}da\notag\\
=&\int_{\partial M}[\langle i_N\omega,J^*(i_{\nabla^{\partial M}V}\omega)\rangle-V\langle i_N\omega, \delta^{\partial M}_f(J^*\omega)\rangle]\,e^{-f}da
\end{align}
from
$\langle J^*\omega,d^{\partial M}V\wedge i_N\omega\rangle=\langle J^*(i_{\nabla^{\partial M}V}\omega),i_N\omega\rangle$.
Therefore, substituting \eqref{1-prop-24} and \eqref{1-prop-25} into \eqref{1-prop-23} yields the desired formula \eqref{1-Th-Formula-1}.

Furthermore, for $\omega\in A^r$ and $S(X_i)=\eta_iX_i$,
$$\aligned
(S^{[p]}\omega)(X_1,\ldots,X_p)=&\sum_{k=1}^p\omega(X_1,\ldots,S(X_k),\ldots,X_p)\\
=&(\eta_1+\ldots+\eta_p)\omega(X_1,\ldots,X_p)
\endaligned$$
gives $S^{[p]}\omega=(\eta_1+\ldots+\eta_p)\omega$ and then
$$*{}_{\partial M}S^{[p]}=(\eta_1+\ldots+\eta_p)*{}_{\partial M}.$$
Similarly, because of $*{}_{\partial M}\omega\in A^{n-p}$, we have
$$
S^{[n-p]}*{}_{\partial M}=(\eta_{p+1}+\ldots+\eta_n)*{}_{\partial M}.
$$
Therefore,
\begin{equation}\label{3-Sec-Savo-9}
*{}_{\partial M}S^{[p]}+S^{[n-p]}*{}_{\partial M}=nH*{}_{\partial M}.
\end{equation}
Since $J^*(*\omega)$ is equal(up to sign) to $*{}_{\partial M}(i_N\omega)$, it follows that:
\begin{equation}\label{3-Sec-Savo-999}\aligned
\langle S^{[n+1-p]}(J^* *\omega),J^* *\omega\rangle=&\langle S^{[n+1-p]}(*{}_{\partial M}(i_N\omega)),*{}_{\partial M}(i_N\omega)\rangle\\
=&\langle[nH*{}_{\partial M}-*{}_{\partial M}S^{[p-1]}](i_N\omega),*{}_{\partial M}(i_N\omega)\rangle\\
=&nH|i_N\omega|^2-\langle S^{[p-1]}(i_N\omega),i_N\omega\rangle
\endaligned\end{equation}
and then $\mathcal{B}_f(\omega,\omega)$ given by \eqref{1-Th-Formula-3} can be written as
\begin{equation}\label{3-Sec-Savo-9999}\aligned
\mathcal{B}_f(\omega,\omega)=&\langle S^{[p]}(J^*\omega),J^*\omega\rangle+nH_f|i_N\omega|^2-\langle S^{[p-1]}(i _N\omega),i_N\omega\rangle\\
=&\langle S^{[p]}(J^*\omega),J^*\omega\rangle+\langle S^{[n+1-p]}(J^* *\omega),J^* *\omega\rangle+f_N|i_N\omega|^2.
\endaligned\end{equation}
This shows that \eqref{1-Th-Formula-3} is equivalent to \eqref{1-Th-Formula-4}

Finally, we have completed the proof of Proposition \ref{prop-2}.

\section{Proof of Theorems \ref{1-Th-2}-\ref{1-Th-6} on compact manifolds}
\subsection{Proof of Theorem \ref{1-Th-2}}

Let $\phi$ be a co-exact $(p-1)$-eigenform associated with $\lambda:=\lambda''_{1,p-1}(\partial M)$. Then from the formula \eqref{dual-1}, the exact $p$-eigenform $\omega=d^{\partial M}\phi$ is also associated with $\lambda$. We consider  $\widehat{\omega}$ satisfying
$$\begin{cases}
d\widehat{\omega}=\delta_f\widehat{\omega}=0,\ \  {\rm on}\ M;\\
J^*\widehat{\omega}=\omega,\ \ \ \ \ \ \ \ \ {\rm on}\ \partial M.
\end{cases}$$
Applying the Reilly formula \eqref{1-Th-Formula-1} with $V=1$ to $\widehat{\omega}$, we obtain
\begin{align}\label{2-thm-1}
0=&\int_M(|\nabla\widehat{\omega}|^2+\langle {W^{[p]}_{f}}(\widehat{\omega}),\widehat{\omega}\rangle)\,e^{-f}dv\notag\\
&+\int_{\partial M}[2\langle\delta^{\partial M}_f(J^*\widehat{\omega}),i_N\widehat{\omega}\rangle+\mathcal{B}_f(\widehat{\omega},\widehat{\omega})]\,e^{-f}da\notag\\
\geq&\int_{\partial M}[2\langle\delta^{\partial M}_f(J^*\widehat{\omega}),i_N\widehat{\omega}\rangle+\mathcal{B}_f(\widehat{\omega},\widehat{\omega})]\,e^{-f}da.
\end{align}
By virtue of \eqref{2-thm-formula-1}, we conclude that
$$\mathcal{B}_f(\widehat{\omega},\widehat{\omega})\geq\sigma_p(\partial M)|J^*\widehat{\omega}|^2+\sigma_{n-p+1}(\partial M)|J^**\widehat{\omega}|^2+f_N|i_N\omega|^2.$$
Since $|J^**\widehat{\omega}|^2=|i_N\widehat{\omega}|^2$ and $\omega=d^{\partial M}\phi$ is an eigenform associated with $\lambda$, that is, $\delta^{\partial M}_f\omega=\delta^{\partial M}_fd^{\partial M}\phi=\lambda\phi$, one has
$$\int_{\partial M}|\omega|^2\,e^{-f}da=\int_{\partial M}|d^{\partial M}\phi|^2\,e^{-f}da=\lambda\int_{\partial M}|\phi|^2\,e^{-f}da$$
and \eqref{2-thm-1} yields
\begin{align}\label{2-thm-2}
0\geq&\int_{\partial M}[2\langle\delta^{\partial M}_f(J^*\widehat{\omega}),i_N\widehat{\omega}\rangle+\mathcal{B}_f(\widehat{\omega},\widehat{\omega})]\,e^{-f}da\notag\\
\geq&\int_{\partial M}[2\lambda \langle\phi,i_N\widehat{\omega}\rangle+\big(\sigma_{n-p+1}(\partial M)+f_N\big)|i_N\widehat{\omega}|^2]\,e^{-f}da\notag\\
&+\sigma_p(\partial M)\int_{\partial M}|\omega|^2\,e^{-f}da.
\end{align}
Substituting the inequality
$$\aligned
2\lambda \langle\phi,i_N\widehat{\omega}\rangle+\big(\sigma_{n-p+1}(\partial M)+f_N\big)|i_N\widehat{\omega}|^2\geq
-\frac{\lambda^2}{\sigma_{n-p+1}(\partial M)+f_N}|\phi|^2
\endaligned$$
into \eqref{2-thm-2} yields
\begin{align}\label{2-thm-222}
0\geq&-\int_{\partial M}\frac{\lambda^2}{\sigma_{n-p+1}(\partial M)+f_N}|\phi|^2\,e^{-f}da+\sigma_p(\partial M)\int_{\partial M}|\omega|^2\,e^{-f}da\notag\\
\geq&-\frac{\lambda^2}{\sigma_{n-p+1}(\partial M)+\inf(f_N)}\int_{\partial M}|\phi|^2\,e^{-f}da+\sigma_p(\partial M)\int_{\partial M}|\omega|^2\,e^{-f}da\notag\\
=&\lambda\Big[\sigma_p(\partial M)-\frac{\lambda}{\sigma_{n-p+1}(\partial M)+\inf(f_N)}\Big]\int_{\partial M}|\phi|^2\,e^{-f}da,
\end{align}
and the desired formula \eqref{2-Th2-1} follows.

\subsection{Proof of Theorem \ref{1-Th-3}}

Let $\omega$ be a $p$-form solution of the absolute cohomology class $H^p(M,\mathbb{R})$. Then $\widehat{\omega}$ satisfies \eqref{3-Th-Formula-1} and the Reilly formula \eqref{1-Th-Formula-1} becomes
\begin{align}\label{3-thm-1}
0=&\int_MV[|\nabla\widehat{\omega}|^2+\langle W^{[p]}_{f,V}(\widehat{\omega}),\widehat{\omega}\rangle]\,e^{-f}dv\notag\\
&+\int_{\partial M}[-V_N|J^*\widehat{\omega}|^2+V\mathcal{B}_f(\widehat{\omega},\widehat{\omega})]\,e^{-f}da.
\end{align}
On $\partial M$, the formula \eqref{1-Th-Formula-3} shows
\begin{align*}
\mathcal{B}_f(\widehat{\omega},\widehat{\omega})=&\langle S^{[p]}(J^*\widehat{\omega}),J^*\widehat{\omega}\rangle\notag\\
\geq&\sigma_p(\partial M)|J^*\widehat{\omega}|^2
\end{align*}
and then \eqref{3-thm-1} gives
\begin{align}\label{3333-thm-1}
0\geq&\int_MV[|\nabla\widehat{\omega}|^2+\langle W^{[p]}_{f,V}(\widehat{\omega}),\widehat{\omega}\rangle]\,e^{-f}dv\notag\\
&+\int_{\partial M}V[\sigma_p(\partial M)-(\ln V)_N]|J^*\widehat{\omega}|^2\,e^{-f}da\geq0.
\end{align}
In particular, if $\sigma_p(\partial M)>(\ln V)_N$, then $\omega$ must be identically zero because it is parallel. Hence, this completes the proof of part $(i)$ on $M$.

For part $(ii)$, there exists a non-trivial solution $\omega$ by assumption, which implies that it must be parallel. Therefore, from \eqref{3333-thm-1} we see that $\sigma_p(\partial M)=(\ln V)_N$.

\subsection{Proof of Theorem \ref{1-Th-5}}
Firstly, we prove the following Pohozhaev-type identity:

\begin{lem}\label{Lem-1}
Let $(M, g, d\mu)$ be a smooth metric measure space. Consider a Lipschitz vector field $F$ on $M$, and  a differential form $\omega\in A^p(M)$. Then the following holds:
\begin{align}\label{1-Lem-1}
\int_M|d\omega|^2{\rm div}_f(F)\,e^{-f}dv=&2\int_M(\langle\nabla F(d\omega),d\omega\rangle-\langle i_Fd\omega,\delta_fd\omega \rangle)\,e^{-f}dv\notag\\
&+\int_{\partial M}(-|d\omega|^2\langle F,N \rangle+2\langle J^*i_Fd\omega,i_Nd\omega \rangle)\,e^{-f}da.
\end{align}
\end{lem}

\proof
It is easy to see
\begin{align}\label{1-lem3.1-1}
{\rm div}_f(|d\omega|^2F)=&\langle F,\nabla|d\omega|^2\rangle+|d\omega|^2{\rm div}_f(F)\notag\\
=&2\langle\nabla_Fd\omega,d\omega\rangle+|d\omega|^2{\rm div}_f(F).
\end{align}
By using the formula \eqref{1-prop-5}, we have
$$\langle d(i_Fd\omega),d\omega\rangle=\langle \nabla_Fd\omega,d\omega\rangle+\langle\nabla F(d\omega),d\omega\rangle$$
and then \eqref{1-lem3.1-1}  becomes
$$|d\omega|^2{\rm div}_f(F)={\rm div}_f(|d\omega|^2F)-2\langle d(i_Fd\omega),d\omega\rangle+2\langle\nabla F(d\omega),d\omega\rangle,$$
which gives
\begin{align*}
\int_M&|d\omega|^2{\rm div}_f(F)\,e^{-f}dv\notag\\
=&\int_M[{\rm div}_f(|d\omega|^2F)-2\langle d(i_Fd\omega),d\omega\rangle+2\langle\nabla F(d\omega),d\omega\rangle]\,e^{-f}dv\notag\\
=&\int_{\partial M}[-|d\omega|^2\langle F,N\rangle+2\langle J^*(i_Fd\omega),i_Nd\omega\rangle]\,e^{-f}da\notag\\
&+2\int_M[\langle\nabla F(d\omega),d\omega\rangle-\langle i_Fd\omega,\delta_fd\omega\rangle]\,e^{-f}dv,
\end{align*}
where we used
$$
\int_M{\rm div}_f(|d\omega|^2F)\,e^{-f}dv=-\int_{\partial M}|d\omega|^2\langle F,N\rangle\,e^{-f}da
$$
and
$$
\int_M\langle d(i_Fd\omega), d\omega\rangle\,e^{-f}dv=\int_M\langle i_Fd\omega,\delta_fd\omega\rangle\,e^{-f}dv-\int_{\partial M}\langle J^*(i_Fd\omega),i_Nd\omega\rangle\,e^{-f}da.
$$
Thus, we finish the proof of Lemma \ref{Lem-1}.
\endproof

Note that Proposition 18 in \cite{Xiong2024} (or see \cite{XX2024}), which states that, for any $\varepsilon>0$, there exists a smooth non-negative function $V_\varepsilon$ such that
$$\nabla^2(-V_\varepsilon)\geq(c-\varepsilon)g$$
on $M\backslash \mathcal{C}$, where $\mathcal{C}$ is an arbitrary fixed neighborhood of cut($\partial M$). Moreover, $V_\varepsilon$ converges uniformly to $V$ as $\varepsilon\rightarrow 0$. Therefore, we can choose $V_\varepsilon$ such that
$$V_\varepsilon=0,\ \nabla V_\varepsilon=N, \ \ \ \ \ \ {\rm on}\  \partial M.$$
Let $\omega$ be an eigenform corresponding to the eigenvalue $\sigma(\partial M)$. By choosing $\omega=d\widehat{\omega}$ and V=$V_\varepsilon$ in \eqref{1-Th-Formula-1},
we get
\begin{align}\label{5-thm-1}
\int_{\partial M}|J^*d\widehat{\omega}|^2\,e^{-f}da=&\int_M[\langle d\widehat{\omega},\nabla^2V_{\varepsilon}(d\widehat{\omega})\rangle +(\Delta_fV_{\varepsilon})|d\widehat{\omega}|^2+V_{\varepsilon}(\langle W^{[p+1]}_f(d\widehat{\omega}),d\widehat{\omega}\rangle\notag\\
&+|\nabla d\widehat{\omega}|^2)]\,e^{-f}dv.
\end{align}
On the other hand, by the Pohozhaev-type identity $\eqref{1-Lem-1}$ with $F=\nabla V_{\varepsilon}$ (because of $F=N$ on $\partial M$ and $\delta_fd\widehat{\omega}=0$ on $M$), we see
\begin{align}\label{5-thm-2}
\int_M&[|d\widehat{\omega}|^2\Delta_fV_{\varepsilon}+2\langle\nabla^2V_{\varepsilon}(d\widehat{\omega}),d\widehat{\omega}\rangle]\,e^{-f}dv\notag\\ =&\int_{\partial M}(|d\widehat{\omega}|^2-2\langle J^*(i_Nd\widehat{\omega}),i_Nd\widehat{\omega}\rangle)\,e^{-f}da\notag\\
=&\int_{\partial M}(|J^*d\widehat{\omega}|^2-|i_Nd\widehat{\omega}|^2)\,e^{-f}da.
\end{align}
Substituting \eqref{5-thm-1} into \eqref{5-thm-2} yields
\begin{align*}
\int_{\partial M}|i_Nd\widehat{\omega}|^2\,e^{-f}da&=\int_M[-\langle\nabla^2V_{\varepsilon}(d\widehat{\omega}),d\widehat{\omega}\rangle
+V_{\varepsilon}(\langle W^{[p+1]}_f(d\widehat{\omega}),d\widehat{\omega}\rangle+|\nabla d\widehat{\omega}|^2)]\,e^{-f}dv\notag\\
&\geq (p+1)(c-\varepsilon)\int_M|d\widehat{\omega}|^2\,e^{-f}dv+\int_MV_{\varepsilon}(\langle W^{[p+1]}_f(d\widehat{\omega}),d\widehat{\omega}\rangle+|\nabla d\widehat{\omega}|^2)\,e^{-f}dv,
\end{align*}
which gives
\begin{align}\label{5-thm-3}
\int_{\partial M}|i_Nd\widehat{\omega}|^2\,e^{-f}da\geq &(p+1)c\int_M|d\widehat{\omega}|^2\,e^{-f}dv+\int_MV(\langle W^{[p+1]}_f(d\widehat{\omega}),d\widehat{\omega}\rangle+|\nabla d\widehat{\omega}|^2)\,e^{-f}dv\notag\\
\geq &(p+1)c\int_M|d\widehat{\omega}|^2\,e^{-f}dv
\end{align}
by letting $\varepsilon\rightarrow 0$.
Applying the definition \eqref{5-Th-formula-2}, we have
\begin{equation}\label{5-thm-4}
\int_{\partial M}|i_Nd\widehat{\omega}|^2\,e^{-f}da=\sigma^2(\partial M)\int_{\partial M}|\omega|^2\,e^{-f}da.
\end{equation}
In addition,
\begin{align*}
0=&\int_M\langle\Delta_f^H\widehat{\omega},\widehat{\omega}\rangle\,e^{-f}dv\notag\\
=&\int_M\langle\delta_fd\widehat{\omega},\widehat{\omega}\rangle\,e^{-f}dv\notag\\
=&\int_M|d\widehat{\omega}|^2\,e^{-f}dv+\int_{\partial M}\langle J^*\widehat{\omega},i_Nd\widehat{\omega}\rangle\,e^{-f}da,
\end{align*}
which yields
\begin{equation}\label{5-thm-5}
\int_M|d\widehat{\omega}|^2\,e^{-f}dv=\sigma(\partial M)\int_{\partial M}|\omega|^2\,e^{-f}da.
\end{equation}
Putting \eqref{5-thm-4} and \eqref{5-thm-5} into \eqref{5-thm-3}, we get the desired inequality \eqref{1-Th-5-formula-1}.

\subsection{Proof of Theorem \ref{1-Th-6}}
It has been proved by Xiong in \cite{Xiong2024} that
$$
(\Delta V_\varepsilon)|d\widehat{\omega}|^2+\langle \nabla^2V_\varepsilon(d\widehat{\omega}),d\widehat{\omega}\rangle
\geq(n-p)(c-\varepsilon)|d\widehat{\omega}|^2.
$$
Then,
$$\aligned
(\Delta_f V_\varepsilon)&|d\widehat{\omega}|^2+\langle \nabla^2V_\varepsilon(d\widehat{\omega}),d\widehat{\omega}\rangle\\
=&(\Delta V_\varepsilon)|d\widehat{\omega}|^2+\langle \nabla^2V_\varepsilon(d\widehat{\omega}),d\widehat{\omega}\rangle+\langle\nabla f,\nabla V_\varepsilon\rangle|d\widehat{\omega}|^2\\
\geq&[(n-p)(c-\varepsilon)+\langle\nabla f,\nabla V_\varepsilon\rangle]|d\widehat{\omega}|^2.
\endaligned$$
Letting $\varepsilon\rightarrow0$ and using Proposition 18 in \cite{Xiong2024} again, then \eqref{5-thm-1} becomes
\begin{align}\label{6-thm-1}
\int_{\partial M}|J^*d\widehat{\omega}|^2\,e^{-f}da\notag\geq&\int_{ M}\{[(n-p)c+\langle\nabla f,\nabla V\rangle]|d\widehat{\omega}|^2+V(|\nabla d\widehat{\omega}|^2\notag\\
&+\langle W^{[p+1]}_f(d\widehat{\omega}),d\widehat{\omega}\rangle)\}\,e^{-f}dv\notag\\
\geq&\int_{ M}[(n-p)c+\langle\nabla f,\nabla V\rangle]|d\widehat{\omega}|^2\,e^{-f}dv.
\end{align}

We recall the function $V=V(\rho)$ defined (see the formulas (3.2) and (3.3) in \cite{XX2024}) by
$$V=\rho-\frac{c}{2}\rho^2$$
with $\max_{M}\rho\leq\frac{1}{c}$, where the distance function $\rho$ is smooth away from the cut locus cut($\partial M$) of $\partial M$. It is easy to see that $|\nabla V|\leq 1$ and \eqref{6-thm-1} gives
\begin{align}\label{6666-thm-1}
\int_{\partial M}|J^*d\widehat{\omega}|^2\,e^{-f}da\geq&\int_{ M}[(n-p)c+\langle\nabla f,\nabla V\rangle]|d\widehat{\omega}|^2\,e^{-f}dv\notag\\
\geq&\int_{ M}[(n-p)c-|\nabla f||\nabla V|]|d\widehat{\omega}|^2\,e^{-f}dv\notag\\
\geq&[(n-p)c-\sup|\nabla f|]\int_{ M}|d\widehat{\omega}|^2\,e^{-f}dv.
\end{align}
Now let $\{\omega_i\}^k_{i=1}$ be the co-closed $p$-forms on $\partial M$ corresponding to the eigenvalues $\{\lambda_i^{[p]}\}_{i=1}^k$, which are orthonormal eigenforms of the weighted Hodge Laplacian with respect to the induced metric on $\partial M$. Let $\{\widehat{\omega}_i\}_{i=1}^p$ be a set of linearly independent solutions to the equation \eqref{5-Th-formula-1} such that $J^*\widehat{\omega}_i=\omega_i$. Then for any $a_i\in\mathbb{R},i=1,2,\ldots,k$, not all zero, it follows from inequality \eqref{6666-thm-1} that
\begin{align*}
\sigma_k^{[p]}\leq&\sup\frac{\int_M|d(\sum_{i=1}^ka_i\widehat{\omega}_i)|^2\,e^{-f}dv}{\int_{\partial M}|J^*(\sum_{i=1}^ka_i\widehat{\omega}_i)|^2\,e^{-f}da}\notag\\
\leq&\frac{1}{(n-p)c-\sup|\nabla f|}\sup\frac{\int_{\partial M}|d^{\partial M}(\sum_{i=1}^ka_i\omega_i)|^2\,e^{-f}da}{\int_{\partial M}|\sum_{i=1}^ka_i\omega_i|^2\,e^{-f}da}\notag\\
\leq&\frac{1}{(n-p)c-\sup|\nabla f|}\lambda_k^{[p]},
\end{align*}
where we used the following fact:
\begin{align*}
\int_{\partial M}\Big|d^{\partial M}\Big(\sum_{i=1}^ka_i\omega_i\Big)\Big|^2\,e^{-f}da
=&\int_{\partial M}\Big\langle \sum_{i=1}^ka_i\omega_i,\delta_f^{\partial M}d^{\partial M}\Big(\sum_{i=1}^ka_i\omega_i\Big)\Big\rangle\,e^{-f}da\notag\\
\leq&\lambda_k^{[p]}\int_{\partial M}\Big|\sum_{i=1}^ka_i\omega_i\Big|^2\,e^{-f}da.
\end{align*}
Now we complete the proof of Theorem \ref{1-Th-6} finally.

\section{Proof of Theorem \ref{1-Th-7} on closed submanifolds}

It is easy to see that for any smooth $p$-form $\omega$ and any smooth function $G$ we have
$$W^{[p]}_f(G\omega)=GW^{[p]}_f(\omega),$$
which from that $W^{[p]}_f$ is $C^\infty(M)$-linear.
Then, \eqref{2-Sec-4} yields
$$[\Delta_f^H,G]\omega=[\nabla^*_f\nabla,G]\omega,$$
which in turn implies 
\begin{align}\label{5-Formula-1}
[\Delta_f^H,G]\omega=&\nabla^*_f\nabla(G\omega)-G\nabla^*_f\nabla\omega\notag\\
=&(\Delta_fG)\omega-2\nabla_{\nabla G}\omega
\end{align}
by virtue of
\begin{align*}
\nabla^*_f\nabla(G\omega)=&(\nabla^*\nabla+\nabla_{\nabla f})(G\omega)\notag\\
=&(\Delta G)\omega-2\nabla_{\nabla G}\omega+G\nabla^*\nabla\omega+\nabla_{\nabla f}(G\omega)\notag\\
=&(\Delta_fG)\omega-2\nabla_{\nabla G}\omega+G\nabla^*_f\nabla\omega.
\end{align*}
Let $G=X_{l}$, where $X_{l}$ is one of the components $X=(X_1,\ldots,X_{M})$, then applying the formula \eqref{5-Formula-1}, we have
\begin{align*}
[[\Delta_f^H,X_{l}],X_{l}]\omega_j=&[\Delta_f^H,X_{l}](X_{l}\omega_j)-X_{l}([\Delta_f^H,X_{l}]\omega_j)\notag\\
=&(\Delta_fX_{l})(X_{l}\omega_j)-2\nabla_{\nabla X_{l}}(X_{l}\omega_j)-X_{l}[(\Delta_fX_{l})\omega_j-2\nabla_{\nabla X_{l}}\omega_j]\notag\\
=&-2|\nabla X_{l}|^2\omega_j,
\end{align*}
and then
$$\int_M|\nabla X_{l}|^2|\omega_j|^2\,e^{-f}dv=-\frac{1}{2}\int_M\langle[[\Delta_f^H,X_{l}],X_{l}]\omega_j,\omega_j\rangle\,e^{-f}dv.$$
Since $\{\lambda_j\}_{j=1}^\infty$ and $\{\omega_j\}_{j=1}^\infty$ are the eigenvalues and orthonormal eigenforms of $\Delta_f^H$ respectively,
using the Theorem 2.2 in \cite{LP2002} (or see Lemma 2.1 in \cite{Ilias2012}) we obtain
\begin{equation}\label{5-Formula-2}
\int_M|\nabla X_{l}|^2|\omega_j|^2\,e^{-f}dv
=\sum_k\frac{(\int_M\langle[\Delta_f^H,X_{l}]\omega_j,\omega_k\rangle\,e^{-f}dv)^2}{\lambda_k^{(p)}-\lambda_j^{(p)}}.
\end{equation}

It is well-known that there exists an orthogonal coordinate system such that the matrix $B=(B_{kl})$,
defined by
$$B_{kl}=\int_M\langle[\Delta_f^H,X_{l}]\omega_j,\omega_{j+k}\rangle\,e^{-f}dv$$
is lower triangular;  this can be achieved via Gram-Schmidt orthogonalization. It follows that $B_{kl}=0$ for $k<l$.

Equation \eqref{5-Formula-2} can be written as follows
\begin{align}\label{5-Formula-3}
\int_M|\nabla X_{l}|^2|\omega_j|^2\,e^{-f}dv=&\sum_{k=1}^{j-1}\frac{(\int_M\langle[\Delta_f^H,X_{l}]\omega_j,\omega_k\rangle\,e^{-f}dv)^2}
{\lambda_k^{(p)}-\lambda_j^{(p)}}\notag\\
&+\sum_{k=j+1}^{j+l-1}\frac{(\int_M\langle[\Delta_f^H,X_{l}]\omega_j,\omega_k\rangle\,e^{-f}dv)^2}
{\lambda_k^{(p)}-\lambda_j^{(p)}}\notag\\
&+\sum_{k=j+l}^{\infty}\frac{(\int_M\langle[\Delta_f^H,X_{l}]\omega_j,\omega_k\rangle\,e^{-f}dv)^2}
{\lambda_k^{(p)}-\lambda_j^{(p)}}.
\end{align}
The first term of the right-hand side of \eqref{5-Formula-3} is non-positive, and the second term is equal to zero via the lower triangle matrix $(B_{kl})$.
Therefore,
\begin{align}\label{5-Formula-4}
\int_M|\nabla X_{l}|^2|\omega_j|^2\,e^{-f}dv\leq&\sum_{k=j+l}^{\infty}\frac{(\int_M\langle[\Delta_f^H,X_{l}]\omega_j,\omega_k\rangle\,e^{-f}dv)^2}
{\lambda_k^{(p)}-\lambda_j^{(p)}}\notag\\
\leq&\frac{1}{\lambda_{j+l}^{(p)}-\lambda_{j}^{(p)}}\sum_{k=j+l}^{\infty}\Big(\int_M\langle[\Delta_f^H,X_{l}]\omega_j,
\omega_k\rangle\,e^{-f}dv\Big)^2\notag\\
\leq&\frac{1}{\lambda_{j+l}^{(p)}-\lambda_{j}^{(p)}}\sum_{k=1}^{\infty}\Big(\int_M\langle[\Delta_f^H,X_{l}]\omega_j,
\omega_k\rangle\,e^{-f}dv\Big)^2.
\end{align}
Summing \eqref{5-Formula-4} over $l$ gives the result:
\begin{align}\label{5-Formula-6}
\sum_{l=1}^{M}(\lambda_{j+l}^{(p)}-\lambda_{j}^{(p)})\Big(\int_M|\nabla X_{l}|^2|\omega_j|^2\,e^{-f}dv\Big)\leq\sum_{l=1}^{M}\int_M|[\Delta_f^H,X_{l}]\omega_j|^2\,e^{-f}dv.
\end{align}
Since X is an isometric immersion, we have $\sum_{l=1}^{M}|\nabla X_{l}|^2=m$ and therefore 
\begin{align}\label{5-Formula-7}
\sum_{l=1}^{M}(\lambda_{j+l}^{(p)}-\lambda_{j}^{(p)})\Big(\int_M|\nabla X_{l}|^2|\omega_j|^2\,e^{-f}dv\Big)=-m\lambda_{j}^{(p)}+\sum_{l=1}^{M}\lambda_{j+l}^{(p)}\int_M|\nabla X_{l}|^2|\omega_j|^2\,e^{-f}dv.
\end{align}
We also have
$$\Delta X=-m\mathbf{H}$$
and
$$
\Delta_fX=-m\mathbf{H}+\nabla f.
$$
It follows that
$$
\sum_{l=1}^{M}(\Delta_fX_{l})^2=m^2|\mathbf{H}|^2+|\nabla f|^2,
$$
$$
\sum_{l=1}^{M}|\nabla_{\nabla X_{l}}\omega_j|^2=|\nabla \omega_j|^2,
$$
$$
\sum_{l=1}^{M}\langle(\Delta_f^HX_{l})\omega_j,\nabla_{\nabla X_{l}}\omega_j\rangle=\frac{1}{2}\langle\nabla  f,\nabla|\omega_j|^2\rangle
$$
and
\begin{align}\label{5-Formula-8}
\sum_{l=1}^{M}&\int_M|[\Delta_f^H,X_{l}]\omega_j|^2\,e^{-f}dv\notag\\
=&\sum_{l=1}^{M}\int_M|(\Delta_fX_{l})\omega_j-2\nabla_{\nabla X_{l}}\omega_j|^2\,e^{-f}dv\notag\\
=&\sum_{l=1}^{M}\int_M[(\Delta_fX_{l})^2|\omega_j|^2+4|\nabla_{\nabla X_{l}}\omega_j|^2-4\langle(\Delta_fX_{l})\omega_j,\nabla_{\nabla X_{l}}\omega_j\rangle]\,e^{-f}dv\notag\\
=&\int_M[(m^2|\mathbf{H}|^2+|\nabla f|^2)|\omega_j|^2+4|\nabla\omega_j|^2-2\langle\nabla f,\nabla|\omega_j|^2\rangle]\,e^{-f}dv\notag\\
\leq&\int_M[(m^2|\mathbf{H}|^2+|\nabla f|^2)|\omega_j|^2+4|\nabla\omega_j|^2]\,e^{-f}dv+2\sup\Delta_ff.
\end{align}
Using \eqref{1-lem-1} again yields
$$\aligned
\int_M|\nabla\omega_j|^2\,e^{-f}dv=&\int_M\langle\nabla_f^*\nabla\omega_j,\omega_j\rangle\,e^{-f}dv\\
=&\int_M[\langle\Delta_f^H\omega_j,\omega_j\rangle-\langle W^{[p]}_f(\omega_j),\omega_j\rangle]\,e^{-f}dv\\
=&\lambda_j^{(p)}-\int_M\langle W^{[p]}_f(\omega_j),\omega_j\rangle\,e^{-f}dv,
\endaligned$$
which transforms \eqref{5-Formula-8} into
\begin{align}\label{5-Formula-9}
&\sum_{l=1}^{M}\int_M|[\Delta_f^H,X_{l}]\omega_j|^2\,e^{-f}dv\notag\\
\leq&4\lambda_j^{(p)}+2\sup\Delta_ff+\int_M[(m^2|\mathbf{H}|^2+|\nabla f|^2)|\omega_j|^2-4\langle W^{[p]}_f(\omega_j),\omega_j\rangle]\,e^{-f}dv.
\end{align}
Now, substituting \eqref{5-Formula-7} and \eqref{5-Formula-9} into \eqref{5-Formula-6}, we obtain
\begin{align}\label{5-Formula-10}
\sum_{l=1}^{M}\lambda_{j+l}^{(p)}\int_M|\nabla X_{l}|^2|\omega_j|^2\,e^{-f}dv\leq&(4+m)\lambda_j^{(p)}+2\sup|\Delta_ff|\notag\\
&+\int_M[(m^2|\mathbf{H}|^2+|\nabla f|^2)|\omega_j|^2-4\langle W^{[p]}_f(\omega_j),\omega_j\rangle]\,e^{-f}dv.
\end{align}
Combining \eqref{5-Formula-10} with the following formula (19) in \cite{Ilias2012}:
$$\sum_{l=1}^{M}\lambda_{j+l}^{(p)}|\nabla X_{l}|^2\geq\sum_{l=1}^{m}\lambda_{j+l}^{(p)},$$
we conclude the proof of Theorem \ref{1-Th-7}.

\section{Appendix}
In this section, we will provide some lemmas on a smooth metric measure space.

\begin{lem}\label{1-App-Lema-1}
Let $(M, g, d\mu)$ be a compact smooth metric measure space with boundary. Then, for $\omega\in A^r(M)$, $\psi\in A^{r+1}(M)$, we have
\begin{equation}\label{3-Appx-1}
\int_M\langle d\omega,\psi\rangle\,e^{-f}dv=\int_M\langle\omega, \delta_f\psi\rangle\,e^{-f}dv-\int_{\partial M}\langle J^*\omega, i_N\psi \rangle\,e^{-f}da.
\end{equation}
\end{lem}

\proof Since the operator $*$ is $C^{\infty}$-linear, we have
$$\aligned
d[\omega\wedge*(e^{-f}\psi)]=&d\omega\wedge*(e^{-f}\psi)-\omega\wedge*[\delta(e^{-f}\psi)]\\
=&e^{-f}d\omega\wedge*\psi-\omega\wedge*[\delta(e^{-f}\psi)],
\endaligned$$
 which is equivalent to
\begin{align}\label{3-Appx-2}
d[\omega\wedge*(e^{-f}\psi)]=&e^{-f}d\omega\wedge*\psi-e^{-f}\omega\wedge*[e^{f}\delta(e^{-f}\psi)]\notag\\
=&e^{-f}d\omega\wedge*\psi-e^{-f}\omega\wedge*(\delta_f\psi).
\end{align}
It follows that
\begin{align}\label{3-Appx-3}
\int_Md[\omega\wedge*(e^{-f}\psi)]=&\int_Me^{-f}d\omega\wedge*\psi-\int_Me^{-f}\omega\wedge*(\delta_f\psi)\notag\\
=&\int_M\langle d\omega,\psi\rangle\,e^{-f}dv-\int_M\langle\omega, \delta_f\psi\rangle\,e^{-f}dv.
\end{align}
By using the Stokes formula, we know that
$$\int_Md[\omega\wedge*(e^{-f}\psi)]=-\int_{\partial M}\langle J^*\omega, i_N(e^{-f}\psi)\rangle\,da,$$
and then \eqref{3-Appx-3} becomes
\begin{align}\label{3-Appx-4}
-\int_{\partial M}\langle J^*\omega, i_N\psi\rangle\,e^{-f}da=\int_M\langle d\omega,\psi\rangle\,e^{-f}dv-\int_M\langle\omega, \delta_f\psi\rangle\,e^{-f}dv,
\end{align}
which is equivalent to \eqref{3-Appx-1}.

\begin{lem}\label{1-App-Lema-2}
Let $(M, g, d\mu)$ be a compact smooth metric measure space. Then, for $\omega\in A^p(M)$, we have
\begin{align}\label{3-Appx2-1}
\int_M(|d\omega|^2+|\delta_f\omega|^2)\,e^{-f}dv=&\int_M\langle\Delta_f^H\omega,\omega\rangle\,e^{-f}dv\notag\\
&+\int_{\partial M}[\langle i_N\omega, J^*(\delta_f\omega)\rangle-\langle J^*\omega, i_Nd\omega\rangle]\,e^{-f}da,
\end{align}
and on $\partial M$, it holds that
\begin{equation}\label{3-Appx2-2}
\delta_f^{\partial M}(J^*\omega)=J^*(\delta_f\omega)+i_N(\nabla_{N}\omega)+S^{[p-1]}(i_N\omega)-nH_f(i_N\omega);
\end{equation}
\begin{equation}\label{3-Appx2-3}
d^{\partial M}(i_{N}\omega)=-i_Nd\omega+J^*(\nabla_{N}\omega)-S^{[p]}(J^*\omega).
\end{equation}
\end{lem}

\proof
By virtue of \eqref{3-Appx-1}, we have
\begin{equation}\label{3-Appx2-4}
\int_M\langle \omega,d\delta_f\omega\rangle\,e^{-f}dv=\int_M\langle \delta_f\omega, \delta_f\omega\rangle\,e^{-f}dv
-\int_{\partial M}\langle J^*\delta_f\omega, i_N\omega \rangle\,e^{-f}da
\end{equation}
and
\begin{equation}\label{3-Appx2-5}
\int_M\langle \omega,\delta_f d\omega\rangle\,e^{-f}dv=\int_M\langle d\omega, d\omega\rangle\,e^{-f}dv+\int_{\partial M}\langle J^*\omega, i_Nd\omega \rangle\,e^{-f}da,
\end{equation}
which shows
\begin{align}\label{3-Appx2-6}
\int_M\langle\Delta_f^H\omega, \omega\rangle\,e^{-f}dv=&\int_M[\langle \omega,d\delta_f\omega\rangle+\langle \omega,\delta_f d\omega\rangle\,]e^{-f}dv\notag\\
=&\int_M(|\delta_f\omega|^2+|d\omega|^2)\,e^{-f}dv-\int_{\partial M}[\langle J^*\delta_{f}\omega, i_N\omega \rangle-\langle J^*\omega, i_Nd\omega \rangle]\,e^{-f}da.
\end{align}
In particular, \eqref{3-Appx2-6} is equivalent to \eqref{3-Appx2-1}.

Using the Gauss formula, we have
\begin{equation}\label{3-Appx2-7}
\nabla_XY=\nabla_X^{\partial M}Y+\langle S(X),Y\rangle N,
\end{equation}
where $\nabla_XN=-S(X)$. Therefore,
$$\aligned
&\nabla_X^{\partial M}(J^*\omega)(Y_1,\ldots,Y_p)\\
=&X[(J^*\omega)(Y_1,\ldots,Y_p)]-\sum_{\alpha=1}^p(J^*\omega)(Y_1,\ldots,\nabla_X^{\partial M}Y_\alpha,\ldots,Y_p)\\
=&X[(J^*\omega)(Y_1,\ldots,Y_p)]-\sum_{\alpha=1}^p(J^*\omega)(Y_1,\ldots,\nabla_XY_\alpha-\langle S(X),Y_\alpha\rangle N,\ldots,Y_p)\\
=&[J^*(\nabla_X\omega)](Y_1,\ldots,Y_p)+[S(X)]^*\wedge i_N\omega(Y_1,\ldots,Y_p),
\endaligned$$
which gives
\begin{equation}\label{3-Appx2-8}
\nabla_X^{\partial M}(J^*\omega)=J^*(\nabla_X\omega)+[S(X)]^*\wedge i_N\omega.
\end{equation}

On the other hand, by the following definition,
$$\aligned
{[i_N(\nabla_X\omega)]}(Y_1,\ldots,Y_{p-1})=&(\nabla_X\omega)(N,Y_1,\ldots,Y_{p-1})\\
=&X[\omega(N,Y_1,\ldots,Y_{p-1})]-\omega(\nabla_XN,Y_1,\ldots,Y_{p-1})\\
&-\sum_{\alpha=1}^{p-1}\omega(N,Y_1,\ldots,\nabla_XY_\alpha,\ldots,Y_{p-1}),
\endaligned$$
we have
\begin{align*}
&\nabla_X^{\partial M}(i_N\omega)(Y_1,\ldots,Y_{p-1})\notag\\
=&X[\omega(N,Y_1,\ldots,Y_{p-1})]-\sum_{\alpha=1}^{p-1}\omega(N,Y_1,\ldots,\nabla_XY_\alpha-\langle S(X),Y_\alpha\rangle N,\ldots,Y_{p-1})\notag\\
=&X[\omega(N,Y_1,\ldots,Y_{p-1})]-\sum_{\alpha=1}^{p-1}\omega(N,Y_1,\ldots,\nabla_XY_\alpha,\ldots,Y_{p-1})\notag\\
=&[i_N(\nabla_X\omega)](Y_1,\ldots,Y_{p-1})-\omega(S(X),Y_1,\ldots,Y_{p-1})\notag\\
=&[i_N(\nabla_X\omega)](Y_1,\ldots,Y_{p-1})-[i_{S(X)}J^*\omega](Y_1,\ldots,Y_{p-1}),
\end{align*}
which gives
\begin{equation}\label{3-Appx2-9}
\nabla_X^{\partial M}(i_N\omega)=i_N(\nabla_X\omega)-i_{S(X)}J^*\omega.
\end{equation}

Thus,
\begin{align*}
&\delta^{\partial M}(J^*\omega)(Y_1,\ldots,Y_{p-1})\notag\\
=&-(\nabla_{e_i}^{\partial M}J^*\omega)(e_i,Y_1,\ldots,Y_{p-1})\notag\\
=&-[J^*\nabla_{e_i}\omega+(S(e_i))^*\wedge i_N\omega](e_i,Y_1,\ldots,Y_{p-1})\notag\\
=&-(\nabla_{e_i}\omega)(e_i,Y_1,\ldots,Y_{p-1})-[(S(e_i))^*\wedge i_N\omega](e_i,Y_1,\ldots,Y_{p-1})\notag\\
=&-e_i[\omega(e_i,Y_1,\ldots,Y_{p-1})]+\omega(\nabla_{e_i}e_i,Y_1,\ldots,Y_{p-1})\notag\\
&+\sum_{\alpha=1}^{p-1}\omega(e_i,Y_1,\ldots,\nabla_{e_i}Y_1,\ldots,Y_{p-1})-\langle S(e_i),e_i\rangle(i_N\omega)(Y_1,\ldots,Y_{p-1})\notag\\
&+\sum_{\alpha=1}^{p-1}(-1)^{\alpha+1}\langle S(e_i),Y_\alpha\rangle(i_N\omega)(e_i,Y_1,\ldots,\widehat{Y_\alpha},\ldots,Y_{p-1})\notag\\
=&(\nabla_{N}\omega)(N,Y_1,\ldots,Y_{p-1})+J^*(\delta\omega)(Y_1,\ldots,Y_{p-1})\notag\\
&-\langle S(e_i),e_i\rangle(i_N\omega)(Y_1,\ldots,Y_{p-1})+\sum_{\alpha=1}^{p-1}(-1)^{\alpha+1}\langle S(e_i),Y_\alpha\rangle(i_N\omega)(e_i,Y_1,\dots,\widehat{Y_\alpha},\ldots,Y_{p-1})\notag\\
=&(\nabla_{N}\omega)(N,Y_1,\ldots,Y_{p-1})+J^*(\delta\omega)(Y_1,\ldots,Y_{p-1})-nH(i_N\omega)(Y_1,\ldots,Y_{p-1})\notag\\
&+\sum_{\alpha=1}^{p-1}(-1)^{\alpha+1}(i_N\omega)(Y_1,\ldots,S(Y_\alpha),Y_2,\ldots,Y_{p-1})\notag\\
=&(\nabla_{N}\omega)(N,Y_1,\ldots,Y_{p-1})+J^*(\delta\omega)(Y_1,\ldots,Y_{p-1})-nH(i_N\omega)(Y_1,\ldots,Y_{p-1})\notag\\
&+[S^{[p-1]}(i_N\omega)](Y_1,Y_2,\ldots,Y_{p-1}),
\end{align*}
which shows
\begin{equation}\label{3-Appx2-10}
\delta^{\partial M}(J^*\omega)=i_N(\nabla_{N}\omega)+J^*(\delta\omega)-nH(i_N\omega)+S^{[p-1]}(i_N\omega),
\end{equation}
where we used the fact:
$$\aligned
J^*(\delta\omega)(Y_1,\ldots,Y_{p-1})=&-(\nabla_{e_A}\omega)(e_A,Y_1,\ldots,Y_{p-1})\\
=&-(\nabla_{N}\omega)(N,Y_1,\ldots,Y_{p-1})-e_i[\omega(e_i,Y_1,\ldots,Y_{p-1})]\\
&+\omega(\nabla_{e_i}e_i,Y_1,\ldots,Y_{p-1})+\sum_{\alpha=1}^{p-1}\omega(e_i,Y_1,\ldots,\nabla_{e_i}Y_\alpha,\ldots,Y_{p-1}).
\endaligned$$

Since $\delta_f=\delta +i_{\nabla f}$, then
\begin{align}\label{3-Appx2-11}
&\delta_f^{\partial M}(J^*\omega)(Y_1,\ldots,Y_{p-1})=\delta^{\partial M}(J^*\omega)(Y_1,\ldots,Y_{p-1})+i_{\nabla^{\partial M} f}(J^*\omega)(Y_1,\ldots,Y_{p-1})\notag\\
=&(\nabla_{N}\omega)(N,Y_1,\ldots,Y_{p-1})+J^*(\delta\omega)(Y_1,\ldots,Y_{p-1})-nH(i_N\omega)(Y_1,\ldots,Y_{p-1})\notag\\
&+[S^{[p-1]}(i_N\omega)](Y_1,Y_2,\ldots,Y_{p-1})+(J^*\omega)(\nabla^{\partial M} f,Y_1,\ldots,Y_{p-1}).
\end{align}
Substituting
$$\aligned
J^*(\delta_f\omega)(Y_1,\ldots,Y_{p-1})=&J^*(\delta\omega)(Y_1,\ldots,Y_{p-1})+J^*(i_{\nabla f}\omega)(Y_1,\ldots,Y_{p-1})\\
=&J^*(\delta\omega)(Y_1,\ldots,Y_{p-1})+(J^*\omega)(\nabla^{\partial M} f,Y_1,\ldots,Y_{p-1})\\
&+f_N(J^*\omega)(N,Y_1,\ldots,Y_{p-1}),
\endaligned$$
with
$$\nabla f=\nabla^{\partial M} f+f_NN,$$
into \eqref{3-Appx2-11} yields
\begin{align*}\label{3-Appx2-12}
&\delta_f^{\partial M}(J^*\omega)(Y_1,\ldots,Y_{p-1})\notag\\
=&(\nabla_{N}\omega)(N,Y_1,\ldots,Y_{p-1})+J^*(\delta_f\omega)(Y_1,\ldots,Y_{p-1})\notag\\
&-nH_f(i_N\omega)(Y_1,\ldots,Y_{p-1})+[S^{[p-1]}(i_N\omega)](Y_1,Y_2,\ldots,Y_{p-1}).
\end{align*}
Therefore, we have
\begin{equation}\label{3-Appx2-13}
\delta_f^{\partial M}(J^*\omega)=J^*(\delta_f\omega)+i_N(\nabla_{N}\omega)+S^{[p-1]}(i_N\omega)-nH_f(i_N\omega)
\end{equation}
and the desired \eqref{3-Appx2-2} is achieved.

By using
$$\aligned
&{[e_i^*\wedge i_N(\nabla_{e_i}\omega)]}(Y_1,\ldots,Y_{p})\\
=&\sum_{\alpha=1}^{p-1}(-1)^{\alpha+1}\langle e_i,Y_\alpha\rangle(\nabla_{e_i}\omega)(N,Y_1,\ldots,\widehat{Y_\alpha},\ldots,Y_{p-1})\\
=&\sum_{\alpha=1}^{p-1}(-1)^{\alpha+1}(\nabla_{Y_\alpha}\omega)(N,Y_1,\ldots,\widehat{Y_\alpha},\ldots,Y_{p-1})
\endaligned$$
and
$$\aligned
-&e_i^*\wedge [i_{S(e_i)}J^*\omega](Y_1,\ldots,Y_{p})\\
=&\sum_{\alpha=1}^{p-1}(-1)^{\alpha}\langle e_i,Y_\alpha\rangle[i_{S(e_i)}J^*\omega](Y_1,\ldots,\widehat{Y_\alpha},\ldots,Y_{p})\\
=&\sum_{\alpha=1}^{p-1}(-1)^{\alpha}\langle e_i,Y_\alpha\rangle(J^*\omega)(S(e_i),Y_1,\ldots,\widehat{Y_\alpha},\ldots,Y_{p})\\
=&-[S^{[p]}(J^*\omega)](Y_1,\ldots,Y_{p}),
\endaligned$$
we have
\begin{align}\label{3-Appx2-14}
d^{\partial M}(i_{N}\omega)(Y_1,\ldots,Y_{p})=&[e_i^*\wedge\nabla_{e_i}^{\partial M}( i_N\omega)](Y_1,\ldots,Y_{p})\notag\\
=&\{e_i^*\wedge [i_N(\nabla_{e_i}\omega)-i_{S(e_i)}J^*\omega]\}(Y_1,\ldots,Y_{p})\notag\\
=&\sum_{\alpha=1}^{p}(-1)^{\alpha+1}(\nabla_{Y_\alpha}\omega)(N,Y_1,\ldots,\widehat{Y_\alpha},\ldots,Y_{p})\notag\\
&-[S^{[p]}(J^*\omega)](Y_1,\ldots,Y_{p}),
\end{align}
where in the second equality, we used \eqref{3-Appx2-9}.
Note that on $\partial M$, we have
\begin{align*}\label{3-Appx2-15}
-(i_Nd\omega)(Y_1,\ldots,Y_{p})=&-d\omega (N,Y_1,\ldots,Y_{p})\notag\\
=&-(e_A^*\wedge \nabla_{e_A}\omega)(N,Y_1,\ldots,Y_{p})\notag\\
=&-(N^*\wedge \nabla_{N}\omega)(N,Y_1,\ldots,Y_{p})-(e_i^*\wedge \nabla_{e_i}\omega)(N,Y_1,\ldots,Y_{p})\notag\\
=&-(\nabla_{N}\omega)(Y_1,\ldots,Y_{p})-\langle e_i,N\rangle(\nabla_{e_i}\omega)(Y_1,\ldots,Y_{p})\notag\\
&+\sum_{\alpha=1}^{p}(-1)^{\alpha+1}\langle e_i,Y_\alpha\rangle(\nabla_{e_i}\omega)(N,Y_1,\ldots,\widehat{Y_\alpha},\ldots,Y_{p})\notag\\
=&-(\nabla_{N}\omega)(Y_1,\ldots,Y_{p})+\sum_{\alpha=1}^{p}(-1)^{\alpha+1}(\nabla_{Y_\alpha}\omega)(N,Y_1,\ldots,
\widehat{Y_\alpha},\ldots,Y_{p})\notag\\
=&-(\nabla_{N}\omega)(Y_1,\ldots,Y_{p})+d^{\partial M}(i_{N}\omega)(Y_1,\ldots,Y_{p})\notag\\
&+[S^{[p]}(J^*\omega)](Y_1,\ldots,Y_{p}),
\end{align*}
where in the last equality, we used \eqref{3-Appx2-14}.
Finally we have the formula \eqref{3-Appx2-3} as follows:
\begin{equation}\label{3-Appx2-16}
d^{\partial M}(i_{N}\omega)=-i_Nd\omega+J^*(\nabla_{N}\omega)-S^{[p]}(J^*\omega). \notag
\end{equation}

\begin{lem}\label{1-App-Lema-3}
$\lambda_{1,p}''=\lambda_{1,p+1}',\ \lambda_{1,p}''=\lambda_{1,n-p}'.$
\end{lem}

\proof
The Hodge decomposition theorem(a similar proof, see \cite{Schwarz1995}) shows
\begin{equation}\label{3-Appx3-1}
A^p=H^p\oplus dA^{p-1}\oplus\delta_f A^{p+1},
\end{equation}
where $H^p=\{\omega | \,\Delta_f^H\omega=0\}$.
Note that, for $\alpha\in A^{p}$ and $d\alpha$ is an eigenform of $\Delta_f^H$(that is, $\Delta_f^H d\alpha=\lambda d\alpha$), then
$$\Delta_f^H(\delta_f d\alpha)=\delta_f \Delta_f^H(d\alpha)=\lambda(\delta_f d\alpha),$$
which shows that $\delta_f(d\alpha)$ is also an eigenform of $\Delta_f^H$. Therefore, we have
\begin{equation}\label{3-Appx3-1111}
\lambda_{1,p}''\leq\lambda_{1,p+1}'.
\end{equation}
Conversely, for $\beta\in A^{p+1}$ and $\delta_f\beta$ is an eigenform of $\Delta_f^H$(that is, $\Delta_f^H (\delta_f\beta)=\lambda \delta_f\beta$), then $d\delta_f\beta$ is also an eigenform of $\Delta_f^H$, which shows
\begin{equation}\label{3-Appx3-2222}
\lambda_{1,p+1}'\leq\lambda_{1,p}''.
\end{equation}
Therefore, \eqref{3-Appx3-1111} and \eqref{3-Appx3-2222} yield
$$\lambda_{1,p}''=\lambda_{1,p+1}'.$$

Next, we prove the second formula. For $\alpha\in A^{p+1}$ and $\delta_f\alpha$ is an eigenform of $\Delta_f^H$(that is, $\Delta_f^H (\delta_f \alpha)=\lambda \delta_f \alpha$), then
$$\Delta_f^H[*(\delta_f \alpha)]=*[\Delta_f^H(\delta_f \alpha)]=\lambda[*(\delta_f \alpha)]$$
and then $*(\delta_f \alpha)$ is also an eigenform of $\Delta_f^H$, which shows that
$$\lambda_{1,n-p}'\leq\lambda_{1,p}''.$$
Conversely, using the similar method, we can obtain
$$\lambda_{1,p}''\leq\lambda_{1,n-p}',$$
and then we have $\lambda_{1,p}''=\lambda_{1,n-p}'$ finally.

\vspace{1cm}
{\bf Acknowledgements} {\ The authors would like to thank anonymous referees for valuable suggestions which made the paper more readable.}

\bibliographystyle{Plain}

\begin{thebibliography}{10}

\bibitem{Bakry85}
D. Bakry, M. Emery,
Diffusion hypercontractives,
S\'{e}m. Prob. XIX. Lect. Notes in Math., 1123(1985), 177-206.

\bibitem{Bakry94}
D. Bakry,
L'hypercontractivit\'{e} et son utilisation en th\'{e}orie des semigroupes,
Lect. Notes Math., 1581(1994), 1-114.

\bibitem{CC2008}
D.G. Chen, Q.-M. Cheng,
Extrinsic estimates for eigenvalues of the Laplace operator,
J. Math. Soc. Japan,  60(2008), 325-339.

\bibitem{CHW2010}
Q.-M. Cheng, G.Y. Huang, G.X. Wei,
Estimates for lower order eigenvalues of a clamped plate problem,
Calc. Var. Partial Differential Equations,  38(2010), 409-416.

\bibitem{Colbois2020}
B. Colbois, A. Girouard, A. Hassannezhad,
The Steklov and Laplacian spectra of Riemannian manifolds with boundary,
J. Funct. Anal.,  278(2020), Paper No. 108409, 38 pp.



\bibitem{Escobar1997}
J.F. Escobar,
The geometry of the first non-zero Stekloff eigenvalue,
J. Funct. Anal.,  150(1997), 544-556.

\bibitem{Escobar1999}
J.F. Escobar,
An isoperimetric inequality and the first Steklov eigenvalue,
J. Funct. Anal., 165(1999), 101-116.



\bibitem{HL2014}
G.Y. Huang, H.Z. Li,
Gradient estimates and entropy formulae of porous medium and fast diffusion equations for the Witten Laplacian,
Pacific J. Math.,  268(2014), 47-78.

\bibitem{HMZ2023}
G.Y. Huang, B.Q. Ma, M.F. Zhu,
Colesanti type inequalities for affine connections,
Anal. Math. Phys.,  13(2023), Paper No. 12, 15 pp.

\bibitem{HMZ2024}
G.Y. Huang, B.Q. Ma, M.F. Zhu,
A Reilly type integral formula and its applications,
Differential Geom. Appl.,  94(2024), Paper No. 102136, 20 pp.

\bibitem{IRS2020}
S. Impera, M. Rimoldi, A. Savo,
Index and first Betti number of $f$-minimal hypersurfaces and self-shrinkers,
Rev. Mat. Iberoam.,  36(2020), 817-840.

\bibitem{Ilias2012}
S. Ilias, O. Makhoul,
A generalization of a Levitin and Parnovski universal inequality for eigenvalues,
J. Geom. Anal.,  22(2012), 206-222.

\bibitem{Karpukhin2017}
M.A. Karpukhin,
Bounds between Laplace and Steklov eigenvalues on nonnegatively curved manifolds,
Electron. Res. Announc. Math. Sci.,  24(2017), 100-109.

\bibitem{Karpukhin2019}
M.A. Karpukhin,
The Steklov problem on differential forms,
Canad. J. Math., 71(2019), 417-435.

\bibitem{Kwong2016}
K.K. Kwong,
Some sharp Hodge Laplacian and Steklov eigenvalue estimates for differential forms,
Calc. Var. Partial Differential Equations,  55(2016), Art. 38, 14 pp.

\bibitem{LP2002}
M. Levitin, L. Parnovski,
Commutators, spectral trace identities, and universal estimates for eigenvalues,
J. Funct. Anal., 192(2002), 425-445.

\bibitem{lW2015}
H.Z. Li, Y. Wei,
$f$-minimal surface and manifold with positive $m$-Bakry-\'{E}mery Ricci curvature,
J. Geom. Anal., 25(2015), 421-435.

\bibitem{lX2019}
J.F. Li, C. Xia,
An integral formula and its applications on sub-static manifolds,
J. Differential Geom.,  113(2019), 493-518.

\bibitem{li05}
X.-D. Li,
Liouville theorems for symmetric diffusion operators on complete Riemannian manifolds,
J. Math. Pures Appl., 84(2005), 1361-1995.

\bibitem{lD11}
X.-D. Li,
Perelman's entropy formula for the Witten Laplacian on Riemannian manifolds via Bakry-\'{E}mery Ricci curvature,
Math. Ann., 353(2012), 403-437.

\bibitem{Peter1996}
P. Li,
Lecture notes on geometric analysis, (1996).

\bibitem{Petersen2012}
P. Petersen, 
Demystifying the Weitzenb\"{o}ck curvature operator,
\href{http://www.math.ucla.edu/~petersen/BLWformulas.pdf}{http://www.math.ucla.edu/~petersen/BLWformulas.pdf}.

\bibitem{QX2015}
G.H. Qiu, C. Xia,
A generalization of Reilly's formula and its applications to a new Heintze-Karcher type inequality,
Int. Math. Res. Not. IMRN,  17(2015), 7608-7619.

\bibitem{Savo2011}
S. Raulot, A. Savo,
A Reilly formula and eigenvalue estimates for differential forms,
J. Geom. Anal., 21(2011), 620-640.

\bibitem{Savo2012}
S. Raulot, A. Savo,
On the first eigenvalue of the Dirichlet-to-Neumann operator on forms,
J. Funct. Anal., 262(2012), 889-914.

\bibitem{Schwarz1995}
G. Schwarz,
Hodge decomposition-a method for solving boundary value problems,
Lecture Notes in Mathematics, 1607. Springer-Verlag, Berlin, 1995.

\bibitem{Sun2008}
H.J. Sun, Q.-M. Cheng, H. Yang,
Lower order eigenvalues of Dirichlet Laplacian,
Manuscripta Math.,  125(2008), 139-156.

\bibitem{XX2024}
C. Xia, C.W. Xiong,
Escobar's conjecture on a sharp lower bound for the first nonzero Steklov eigenvalue,
Peking Math. J.,  7(2024), 759-778.

\bibitem{Wang2018}
C.Y. Xia, Q.L. Wang,
Eigenvalues of the Wentzell-Laplace operator and of the fourth order Steklov problems,
J. Differential Equations, 264(2018), 6486-6506.


\bibitem{Xiong2018}
C.W. Xiong,
Comparison of Steklov eigenvalues on a domain and Laplacian eigenvalues on its boundary in Riemannian manifolds,
J. Funct. Anal.,  275(2018), 3245-3258.

\bibitem{Xiong2024}
C.W. Xiong,
A weighted Reilly formula for differential forms and sharp Steklov eigenvalue estimates,
Sci. China Math., (2025),
https://doi.org/10.1007/s11425-024-2432-1.

\bibitem{Wang2009}
Q.L. Wang, C.Y. Xia,
Sharp bounds for the first non-zero Stekloff eigenvalues,
J. Funct. Anal., 257(2009), 2635-2644.

\bibitem{Wei09}
G.-F. Wei, W. Wylie,
Comparison geometry for the Bakry-Emery Ricci tensor,
J. Differ. Geom., 83(2009), 377-405.


\end{thebibliography}

\end{document}